\documentclass[12pt]{article}
\usepackage{latexsym,amssymb,amsmath}
\begin{document}

\baselineskip=20pt

\newcommand{\la}{\langle}
\newcommand{\ra}{\rangle}
\newcommand{\psp}{\vspace{0.4cm}}
\newcommand{\pse}{\vspace{0.2cm}}
\newcommand{\ptl}{\partial}
\newcommand{\dlt}{\delta}
\newcommand{\sgm}{\sigma}
\newcommand{\al}{\alpha}
\newcommand{\be}{\beta}
\newcommand{\G}{\Gamma}
\newcommand{\gm}{\gamma}
\newcommand{\vs}{\varsigma}
\newcommand{\Lmd}{\Lambda}
\newcommand{\lmd}{\lambda}
\newcommand{\td}{\tilde}
\newcommand{\vf}{\varphi}
\newcommand{\yt}{Y^{\nu}}
\newcommand{\wt}{\mbox{wt}\:}
\newcommand{\rd}{\mbox{Res}}
\newcommand{\for}{\mbox{for}}
\newcommand{\ad}{\mbox{\small ad}}
\newcommand{\voa}{(V,Y(\cdot,z),{\bf 1},\omega)}
\newcommand{\stl}{\stackrel}
\newcommand{\ol}{\overline}
\newcommand{\ul}{\underline}
\newcommand{\es}{\epsilon}
\newcommand{\dmd}{\diamond}
\newcommand{\clt}{\clubsuit}
\newcommand{\vt}{\vartheta}
\newcommand{\ves}{\varepsilon}
\newcommand{\dg}{\dagger}
\newcommand{\Dlt}{\Delta}
\newcommand{\Edo}{\mbox{End}}

\begin{center}{\LARGE \bf Simple Conformal Superalgebras}\end{center}
\begin{center}{\LARGE \bf of Finite Growth}
\end{center}
\vspace{0.2cm}

\begin{center}{\large Xiaoping Xu}\end{center}
\begin{center}{Department of Mathematics, The Hong Kong University of Science \& Technology}\end{center}
\begin{center}{Clear Water Bay, Kowloon, Hong Kong}\footnote{Research supported
 by Hong Kong RGC Competitive Earmarked Research Grant HKUST709/96P.}\end{center}

\vspace{0.3cm}

\begin{center}{\Large \bf Abstract}\end{center}
\vspace{0.2cm}

{\small In this paper, we construct six families of infinite simple conformal superalgebra of finite growth based on our earlier work on constructing vertex operator superalgebras from graded assocaitive algebras. Three subfamilies of these conformal superalgebras are generated by simple Jordan algebras of types A, B and C in a certain sense}.

\section{Introduction}

The notion of conformal superalgebra was formulated by Kac [K3]. Conformal superalgebras play important roles in quantum field theory (e.g. cf. [K3]) and vertex operator superalgebras (e.g. cf. [K3], [X5]). The classification theorem of simple conformal superalgebras of finite type was announced by Kac [K5] and proved in [DK]. Except the algebra $CK_6$, all the classified finite simple conformal superalgebras are essentially quite known (e.g. cf. [K3], [K4]). The algebra $CK_6$ is a subalgebra of the algebra $K_6$ constructed through the Hodge dual (cf. [CK1]). A natural question is whether there exist new simple conformal superalgebras whose structures are close to the simple conformal superalgebras of finite type. In this paper, we shall give an affirmative answer.  

Motivated by the vertex operator subalgebras generated by certain quadratic free fields in our earlier work [X3]  on ternary moonshine spaces, we  introduced in Section 7.3 of [X5] a new family of infinite-dimensional Lie superalgebras, which we called ``double affinizations'' of $\Bbb{Z}_2$-graded associative algebras with respect to a trace map. From these Lie superalgebras, we constructed new families of conformal superalgebras with a Virasoro element, which yielded new families of simple vertex operator superalgebras generated by their subspaces of small weights. In this paper, we shall construct six families of infinite simple conformal superalgebras of finite growth from matrix algebras and prove their simplicity. A generator subset of each of these algebras is determined. In particular, three subfamilies of these algebras are simple conformal algebras generated by simple Jordan algebras of types A, B and C in a certain sense except their ``minimal cases." Below we shall give more detailed introduction.

Throughout this chapter, the base field $\Bbb{F}$ is an arbitrary field of characteristic 0. 
 For two vector spaces $V_1$ and $V_2$, we denote by $LM(V_1,V_2)$ the space of linear maps from $V_1$ to $V_2$. Moreover, we denote by $\Bbb{Z}$ the ring of integers, by $\Bbb{N}$ the set of natural numbers $\{0,1,2,....\}$ and by $\Bbb{Z}_2=\Bbb{Z}/2\Bbb{Z}$ the cyclic group of order 2. When the context is clear, we use $\{0,1\}$ to denote the elements of $\Bbb{Z}_2$. We shall also use the following operator of taking residue:
$$\rd_z(z^n)=\dlt_{n,-1}\qquad\for\;\;n\in \Bbb{Z}.\eqno(1.1)$$
Furthermore, all the binomials are assumed to be expanded in the nonnegative powers of the second variable. 

 A {\it conformal superalgebra} $R=R_0\oplus R_1$ is a $\Bbb{Z}_2$-graded $\Bbb{C}[\ptl]$-module with a $\Bbb{Z}_2$-graded linear map $Y^+(\cdot,z):\;R\rightarrow LM(R,R[z^{-1}]z^{-1})$ satisfying:
$$Y^+(\ptl u,z)={dY^+(u,z)\over dz}\qquad\for\;\;u\in R;\eqno(1.2)$$
$$Y^+(u,z)v=(-1)^{ij}\rd_x{e^{x\ptl}Y^+(v,-x)u\over z-x},\eqno(1.3)$$
$$Y^+(u,z_1)Y^+(v,z_2)-(-1)^{ij}Y^+(v,z_2)Y^+(u,z_1)=\rd_x{Y^+(Y^+(u,z_1-x)v,x)\over z_2-x} \eqno(1.4)$$
for $u\in R_i;\;v\in R_j$. We denote by $(R,\ptl,Y^+(\cdot,z))$ a conformal superalgebra. When $R_1=\{0\}$, we simply call $R$ a {\it conformal algebra}. 

The above definition is the equivalent generating-function form to that given in [K3], where the author used the component formulae with $Y^+(u,z)=\sum_{n=0}^{\infty}u_{(n)}z^{-1}$.

An {\it ideal} ${\cal I}$ of a conformal superalgebra $(R,\ptl,Y^+(\cdot,z))$ is a subspace of $R$ such that
$$\ptl({\cal I})\subset {\cal I},\;\;Y^+(u,z)({\cal I})\subset{\cal I}[z^{-1}]\qquad\for\;\;u\in R.\eqno(1.5)$$
By (1.3), we have
$$Y^+({\cal I},z)R\subset {\cal I}[z^{-1}]\eqno(1.6)$$
for an ideal ${\cal I}$ of $R$. The conformal superalgebra $R$ is called {\it simple} if the only ideals of $R$ are $\{0\}$ and $R$.

 Let $\G$ be an additive subgroup of $\Bbb{F}$ such that $\G$ has a subgroup $\G_0$ of index 2. A conformal superalgebra $(R,\ptl,Y^+(\cdot,z))$ is called $\G$-{\it weighted} if for an index-2 subgroup $\G_0$ of $\G$ and its coset $\G_1=\G\setminus\G_0$,
$$R_i=\bigoplus_{\al\in\G_i}R^{(\al)}\;\;\mbox{as subspaces}\qquad\for\;\;i=0,1\eqno(1.7)$$
such that for $u\in R^{(\mu)},\;Y^+(u,z)=\sum_{n\in\Bbb{N}-\mu+1}u(n)z^{-n-\mu}$,
$$\ptl R^{(\al)}\subset R^{(\al-1)},\;\;\;u(m)R^{(\al)}\subset R^{(\al-m)}.\eqno(1.8)$$
The elements in $R^{(\al)}$ are  called  the  {\it elements of weight} $\al$ and $R^{(\al)}$ is called the {\it subspace of weight} $\al$. A $\G$-weighted conformal superalgebra $(R,Y^+(\cdot,z))$ is said of {\it finite growth} if it contains a subspace 
$$V=\bigoplus_{\al\in\G}V^{(\al)},\qquad V^{(\al)}=V\bigcap R^{(\al)}\eqno(1.9)$$ 
such that
$$R=\Bbb{F}[\ptl]V\;\;\mbox{and}\;\;\mbox{dim}\:V^{(\al)}<N_0\;\;\for\;\;\al\in\G,\eqno(1.10)$$
where $N_0$ is a fixed positive integer. For a $\G$-weighted conformal superalgebra $(R,Y^+(\cdot,z))$, we define the {\it weight system}:
$$\Dlt=\{\al\in\G\mid R^{(\al)}\neq\{0\}\}.\eqno(1.11)$$
We shall also called the conformal superalgebra $(R,Y^+(\cdot,z))$ $\Dlt$-{\it weighted} (the weighting group is clearly the additive subgroup of $\Bbb{F}$ generated by $\Dlt$). In nonsuper case ($R_1=\{0\}$), we say $R$ is $\G_0$-{\it weighted} ($\G_1$ is redundant). Since $\Bbb{Z}$ is the only index-2 subgroup of $\Bbb{Z}/2$, we have $\G_0=\Bbb{Z}$ when we consider a $\Bbb{Z}/2$-weighted conformal superalgebra.

A sub-superalgebra $R'$ of a  conformal superalgebra $(R,Y^+(\cdot,z))$ is said to {\it be generated by a subset} $S$ if
$$R'=\mbox{span}\:\{u_{m_1}^1\cdots u^p_{m_p}v\mid u^j,v\in S,\;p,m_j\in\Bbb{N}\},\eqno(1.12)$$
where we write
$$Y^+(u,z)=\sum_{n=0}^{\infty}u_nz^{-n-1}\qquad\for\;\;u\in R.\eqno(1.13)$$

In this paper, we shall construct six families of $(1+\Bbb{N}/2)$-weighted simple conformal superalgebras of finite growth. Three subfamilies of these algebras are generated by their subspaces of minimal weight, whose homogeneous structures are simple Jordan algebras of types A, B and C except their ``minimal cases." 

Our first family of simple conformal algebras $R_{k\times k,\ell}$ are parametrized by two positive integral variables, related to the algebra of $k\times k$ matrices. The algebra $R_{1\times 1,2}$ is the well-known $W_{\infty}$ algebra  without center (cf. [Ba]) and the algebra $R_{1\times 1,1}$ is the well-known $W_{1+\infty}$ algebra without center (cf. [PRS]) in mathematical physics. The more general algebra $R_{k\times k,1}$ is the $W_{1+\infty}(gl_k)$ algebra studied by van de Leur [V] without center related to $k$ component KP hierarchy. Our first family of simple conformal superalgebras $R_{[k_1,k_2],\ell}$ are parametrized by three positive integral variables, related to the algebra of $(k_1+k_2)\times (k_1+k_2)$ matrices. We believe that the algebra $R_{[1,1],1}$ is related to the supersymmetric analogues of the $W_{1+\infty}$ algebra  studied by physicists [DHP], [Y] and [YW]. It is conceivable that all the simple conformal superalgebras presented in this paper would eventually be related to certain integrable systems that are generalizations of KP hierarchy.

In Section 2, we shall present the general construction of conformal superalgebras from $\Bbb{Z}_2$-graded associative algebras and its motivation from quadratic free fields. In Section 3, we shall construct three families of infinite simple conformal algebras of finite growth from matrix algebras and prove their simplicity. Section 4 is devoted to the constructions of three families of infinite simple conformal superalgebras of finite growth with nonzero odd part from  matrix algebras and the proof of their simplicity.

\section{Motivation and General Construction}

In this section, we shall first give a motivation from quadratic free bosonic fields of the construction. Then we present our general construction of conformal superalgebras from $\Bbb{Z}_2$-graded assocaitive algebras based on Section 7.3 of [X5]. Moreover, we shall make a comparison of the general construction with the conformal superalgebras that generate the loop algebras and the centerless Virasoro algebras. 

\subsection{Motivation}

In this subsection, we shall single out the conformal subalgebras related to certain quadratic free fields in our earlier works [X3] on ternary moonshine spaces.

Let $H$ be a vector space with a nondegenerate symmetric bilinear form $\la\cdot,\cdot\ra$ such that there exist two subspaces $H_+,H_-$ satisfying $H=H_++H_-$ and
$$\la H_+,H_+\ra=\la H_-,H_-\ra=\{0\}.\eqno(2.1)$$
Let $t$ be an indeterminate and set
$$\hat{H}=H\otimes_{\Bbb{F}}\Bbb{F}[t,t^{-1}]\oplus \Bbb{F}\kappa,\eqno(2.2)$$
where $\kappa$ is a symbol to denote a base vector of one-dimensional vector space. We define algebraic operation $[\cdot,\cdot]$ on $\hat{H}$ by
$$[h_1\otimes t^m+\lmd_1\kappa,h_2\otimes t^n+\lmd_2\kappa]=m\la h_1,h_2\ra\dlt_{m+n,0}\kappa\eqno(2.3)$$
for $h_1,h_2\in H,\;m,n\in\Bbb{Z},\;\lmd_1,\lmd_2\in\Bbb{F}$. Then $(\hat{H},[\cdot,\cdot])$ forms a Lie algebra, which is called a {\it Heisenberg Lie algebra}.\index{Heisenberg Lie algebra} For convenience, we denote
$$h(m)=h\otimes t^m\qquad\for\;\;h\in H,\;m\in\Bbb{Z}.\eqno(2.4)$$
Set
$$\hat{H}_-=\mbox{span}\:\{h(- m)\mid h\in H,\;m\in\Bbb{Z}_+\},\;\;\hat{B}_H=\mbox{span}\:\{\kappa,h(m)\mid h\in H,\;m\in\Bbb{N}\}.\eqno(2.5)$$
Then $\hat{H}_-$ and $\hat{B}_H$ are trivial Lie subalgebras of $\hat{H}$ and
$$\hat{H}=\hat{H}_-\oplus \hat{B}_H.\eqno(2.6)$$

Let $\Bbb{F}{\bf 1}$ be a one-dimensional vector space with the base element ${\bf 1}$. We define an action of $\hat{B}_H$ on 
$\Bbb{F}{\bf 1}$ by
$$h(m)({\bf 1})=0,\;\;\kappa({\bf 1})={\bf 1}\qquad\for\;\;h\in H,\;m\in\Bbb{N}.\eqno(2.7)$$
Then $\Bbb{F}{\bf 1}$ forms a $\hat{B}_H$-module. We denote by $U(\cdot)$ the universal envelopping algebra of a Lie algebra and by $S(\cdot)$ the symmetric algebra generated by a vector space. Form an induced $\hat{H}$-module
$$V=U(\hat{H})\otimes_{U(\hat{B}_H)}\Bbb{F}{\bf 1}\;\;(\cong S(\hat{H}_-)\otimes_{\Bbb{F}}\Bbb{F}{\bf 1}\;\mbox{as vector spaces}).\eqno(2.8)$$
Moreover, we set
$$h^+(z)=\sum_{m=0}^{\infty}h(m)z^{-m-1},\;\;h^-(z)=\sum_{m=1}^{\infty}h(-m)z^{m-1},\;\;h(z)=h^+(z)+h^-(z)\eqno(2.9)$$
for $h\in H$, where $z$ is a formal variable. As operators on $V$, $\{h(z)\mid h\in H\}$ are called {\it free bosonic fields}.\index{free bosonic fields}

For convenience, we denote
$$u\otimes{\bf 1}=u\qquad\for\;\;u\in S(\hat{H}_-).\eqno(2.10)$$
Set
$$\hat{R}_2=\mbox{span}\:\{h_1(-m_1)h_2(-m_2),{\bf 1}\mid h_1,h_2\in H,\;m_1,m_2\in\Bbb{Z}_+\}.\eqno(2.11)$$
We define a linear map $Y(\cdot,z):\hat{R}_2\rightarrow LM(V,V[z^{-1},z]])$ by
$$Y(h_1(-m-1)h_2(-n-1),z)={1\over m!n!}\left({d^mh^-_1(z)\over dz^m}{d^nh_2(z)\over dz^n}+{d^nh_2(z)\over dz^n}{d^mh^+_1(z)\over dz^m}\right)\eqno(2.12)$$
for $h_1,h_2\in H$ and $m,n\in\Bbb{N}$ and
$$Y({\bf 1},z)=\mbox{Id}_V.\eqno(2.13)$$
The operator $Y(h_1(-m-1)h_2(-n-1),z)$ is called a {\it quadratic bosonic field}.\index{quadratic bosonic field} Moreover, we write
$$Y(u,z)=\sum_{n\in\Bbb{Z}}u_nz^{-n-1},\;\;Y^+(u,z)=\sum_{n=0}^{\infty}u_nz^{-n-1}\qquad\for\;\;u\in\hat{R}_2.\eqno(2.14)$$
In particular, 
\begin{eqnarray*}& &(h_1(-m-1)h_2(-n-1))_k\\&=&\sum_{j=0}^{\infty}(^{-j-1}_{\;\;\;n})[(^{j+m+n-k}_{\;\;\;\;\;\;\;\:m})h_1(k-m-n-j-1)h_2(j)\\& &+(^{-j-1}_{\;\;\;m})(^{j+m+n-k}_{\;\;\;\;\;\;\;n})h_2(k-m-n-j-1)h_1(j)]\hspace{5.7cm}(2.15)\end{eqnarray*}
for $h_1,h_2\in H$ and $m,n,k\in\Bbb{N}$.

 Note that
\begin{eqnarray*}&&h_1(m)h_2(n)(h_3(-j)h_4(-k))\\&=&mn(\dlt_{m,j}\dlt_{n,k}\la h_1,h_3\ra \la h_2,h_4\ra+\dlt_{m,k}\dlt_{n,j}\la h_1,h_4\ra \la h_2,h_3\ra){\bf 1},\hspace{4.2cm}(2.16)\end{eqnarray*}
\begin{eqnarray*}& &h_1(-m)h_2(n)(h_3(-j)h_4(-k))\\&=&n\dlt_{n,j}\la h_2,h_3\ra h_1(-m)h_4(-k)+n\dlt_{n,k}\la h_2,h_4\ra h_1(-m)h_3(-j)\hspace{3.7cm}(2.17)\end{eqnarray*}
for $h_1,h_2,h_3,h_4\in H$ and $m,n,j,k\in\Bbb{Z}_+$. Expressions (2.15)-(2.17) show that
$$Y^+(u,z)v\subset \hat{R}_2[z^{-1}]\qquad\for\;\;u,v\in \hat{R}_2.\eqno(2.18)$$
Moreover, we define $\ptl\in\Edo \hat{R}_2$ by
$$\ptl({\bf 1})=0,\;\;\ptl (h_1(-m)h_2(-n))=mh_1(-m-1)h_2(-n)+nh_1(-m)h_2(-n-1)\eqno(2.19)$$
for $h_1,h_2\in H$ and $m,n\in\Bbb{Z}_+$. Then the family $(\hat{R}_2,\ptl,Y^+(|_{\hat{R}_2},z))$ forms a conformal algebra by (3.3.42) and Theorem 6.1.3 in [X5].

According to linear algebra, there exist basis $\{\vs^{\pm}_j|j\in I\}$ of $H_{\pm}$ (cf. (2.1)) such that
$$\la\vs^+_i,\vs_j^-\ra=\dlt_{i,j}\qquad\for\;\;i,j\in I\eqno(2.20)$$
by (2.1) and nondegeneracy of $\la\cdot,\cdot\ra$, where $I$ is an index set. Note that
$$\vs^+_{j_1}(-1)\vs^-_{j_2}(1)(\vs^+_{j_3}(-1)\vs^-_{j_4}(-1))=\dlt_{j_2,j_3}\vs^+_{j_1}(-1)\vs^-_{j_4}(-1),\eqno(2.21)$$
$$\vs^-_{j_1}(-1)\vs^+_{j_2}(1)(\vs^+_{j_3}(-1)\vs^-_{j_4}(-1))=\dlt_{j_2,j_4}\vs^+_{j_3}(-1)\vs^-_{j_1}(-1)\eqno(2.22)$$
for $j_1,j_1,j_3,j_4\in I$. Expressions (2.21) and (2.22) are essentially equivalent to matrix multiplications! This shows that there exists the connection between matrix algebra of dimension $I\times I$ and the conformal algebra $(\hat{R}_2,\ptl,Y^+(|_{\hat{R}_2},z))$, which is a motivation of our general construction of conformal superalgebra from $\Bbb{Z}_2$-graded associative algebras.

\subsection{General Construction}

In this subsection, we shall present the general construction of conformal superalgebra from $\Bbb{Z}_2$-graded associative algebras. 

Let 
$${\cal A}={\cal A}_0\oplus {\cal A}_1\eqno(2.23)$$
 be a $\Bbb{Z}_2$-graded associative algebra with an identity element $e$.  
Let $M_{2\times 2}({\cal A})$ be the algebra of $2\times 2$ matrices whose entries are in ${\cal A}$. Note that we have the following subalgebra of $M_{2\times 2}({\cal A})$:
$${\cal N}({\cal A})=\left\{\left(\begin{array}{cc}a_{0,0},&a_{0,1}\\ a_{1,0},&a_{1,1}\end{array}\right)\mid a_{0,0},a_{1,1}\in{\cal A}_0,\;a_{0,1},a_{1,0}\in{\cal A}_1\right\}.\eqno(2.24)$$
Set
$$R({\cal A})={\cal N}({\cal A})\otimes_{\Bbb{F}}\Bbb{F}[t_1,t_2],\eqno(2.25)$$
where $t_1$ and $t_2$ are indeterminates. Denote
$$u(n_1,n_2)=u\otimes t^{n_1}_1t_2^{n_2}\;\;\;\for\;\;u\in {\cal N}({\cal A}),\;n_1,n_2\in \Bbb{N}.\eqno(2.26)$$
We make a convention that any notions that appear technically and have not been defined are treated as zero. For instance,
$$v(-1,0)=0,\;\;w(2,-3)=0\qquad\;\;\mbox{if}\;\;v,w\in {\cal N}({\cal A}).\eqno(2.27)$$
For convenience, we denote
$$u_{[0,0]}=\left(\begin{array}{cc}u,&0\\ 0,&0\end{array}\right),\qquad 
u_{[1,1]}=\left(\begin{array}{cc}0,&0\\ 0,&u\end{array}\right)\qquad\for\;\;u\in{\cal A}_0,\eqno(2.28)$$
$$v_{[0,1]}=\left(\begin{array}{cc}0,&v\\ 0,&0\end{array}\right),\qquad
 v_{[1,0]}=\left(\begin{array}{cc}0,&0\\ v,&0\end{array}\right)\qquad\for\;\;v\in{\cal A}_1.\eqno(2.29)$$
Moreover, we make a convention that if the notion $u_{[i,j]}$ is used for $i,j\in\Bbb{Z}_2$, we always mean $u\in {\cal A}_{i+j}$. Under this convention, we have 
$$R({\cal A})=\mbox{span}\:\{u_{[i,j]}(m,n)\mid i,j\in\Bbb{Z}_2,\;u\in{\cal A},\;m,n\in\Bbb{N}\}.\eqno(2.30)$$

Now we define a linear map $Y^+(\cdot,z):R({\cal A})\rightarrow R({\cal A})[z^{-1}]z^{-1}$ by
\begin{eqnarray*}& &Y^+(u_{[i_1,i_2]}(m_1,m_2),z)v_{[j_1,j_2]}(n_1,n_2)\\&=&
\dlt_{i_2,j_1}(n_1\dlt_{i_2,0}+1)(^{-n_1-1-\dlt_{i_2,0}}_{\;\;\;\;\;\;\;m_2})\sum_{p=0}^{m_1+m_2+n_1+\dlt_{i_2,0}}(^{\;p}_{m_1})(uv)_{[i_1,j_2]}(p,n_2)z^{p-m_1-m_2-n_1-1-\dlt_{i_2,0}}\\&&-(-1)^{(i_1+i_2)(j_1+j_2)}\dlt_{i_1,j_2}(\dlt_{i_1,1}-(n_2+1)\dlt_{i_1,0})(^{-n_2-1-\dlt_{i_1,0}}_{\;\;\;\;\;\;\;m_1})\\& &\sum_{q=0}^{m_1+m_2+n_2+\dlt_{i_1,0}}(^{\;q}_{m_2})(vu)_{[j_1,i_2]}(n_1,q)z^{q-m_1-m_2-n_2-1-\dlt_{i_1,0}}\hspace{4.6cm}(2.31)\end{eqnarray*}
for $u_{[i_1,i_2]}(m_1,m_2),v_{[j_1,j_2]}(n_1,n_2)\in R({\cal A})$. The above formula was motivated by the following formula 
of the nonnegative operators of a quadratic field action on a quadratic element:
$$Y^+(\vs_{j_1}^+(-m_1-1)\vs^-_{j_2}(-m_2-1),z)\vs_{j_3}^+(-n_1-1)\vs^-_{j_4}(-n_2-1)\qquad \eqno(2.32)$$
for $j_1,j_2,j_3,j_4\in I$ and $m_1,n_1,m_2,n_2\in\Bbb{N}$, which is defined by (2.1), (2.12), (2.14) and (2.20). The reader may calculate (2.32) by (2.15) and compare it with (2.31) when $i_1=i_2=j_1=j_2=0$ in order to better understand (2.31). The action of $\Bbb{F}[\ptl]$ on $R({\cal A})$ is defined by
$$\ptl u_{[i,j]}(m,n)=(m+1)u_{[i,j]}(m+1,n)+(n+1)u_{[i,j]}(m,n+1)\eqno(2.33)$$
for $i,j\in\Bbb{Z}_2,\;u\in{\cal A}_{i+j}$ and $m,n\in\Bbb{N}$. By Theorem 7.3.4 in [X5], $(R({\cal A}),\ptl,Y^+(\cdot,z))$ forms a $(1+\Bbb{N}/2)$-weighted conformal superalgebra with
$$R({\cal A})^{(n)}=\mbox{span}\:\{u_{[i,j]}(m_1,m_2)\in R({\cal A})\mid m_1+m_2+\dlt_{i,0}+\dlt_{j,0}+(\dlt_{i,1}+\dlt_{j,1})/2=n\}\eqno(2.34)$$
for $n\in\Bbb{N}/2$. The correspondence between the notations in the above and the those in Section 7.3 of [X5] is as follows: the conformal superalgebra $R({\cal A})$ in the above is the quotient algebra $R({\cal A})/\Bbb{F}\ul{1}$ in Section 7.3 of [X5] (cf. (7.3.51)) and the notion
$$ u_{[i,j]}[m,n]=u(-1+i/2-m,-1+j/2-n)\otimes\ul{1}+\Bbb{F}\ul{1}.\eqno(2.35)$$
 Furthermore, $R({\cal A})$ is a free $\Bbb{F}[\ptl]$-module over the subspace
$$V_{\cal A}=\mbox{span}\:\{u_{[i,j]}(0,m)\mid u\in{\cal A},\;i,j\in\Bbb{Z}_2,\;m\in\Bbb{N}\}\eqno(2.36)$$
and
$$\mbox{dim}\;V_{\cal A}\bigcap R({\cal A})^{(n)}\leq\max\{2\mbox{dim}\:{\cal A}_0,\;2\mbox{dim}\:{\cal A}_1\}.\eqno(2.37)$$
Thus $R({\cal A})$ is a $(1+\Bbb{N}/2)$-weighted conformal superalgebra of finite growth if ${\cal A}$ is finite-dimensional. So are its subalgebras.
  
Let
$$u_{[i,j]}(x_1,x_2)=\sum_{m,n=0}^{\infty}u_{[i,j]}(m,n)x_1^mx_2^n\qquad\for\;\;i,j\in\Bbb{Z}_2,\;u\in{\cal A}_{i+j},\eqno(2.38)$$
where $x_1$ and $x_2$ are formal variables. Then (2.30) can be rewritten as
\begin{eqnarray*}&&Y^+(u_{[i_1,i_2]}(x_1,x_2),z)v_{[j_1,j_2]}(y_1,y_2)\\&=&\dlt_{i_2,j_1}\ptl_{y_1}^{\dlt_{i_2,0}}\left[{(uv)_{[i_1,j_2]}(y_1-x_2+x_1,y_2)\over z+x_2-y_1}\right]\\& &+(-1)^{(i_1+i_2)(j_1+j_2)+j_2}\dlt_{j_2,i_1}\ptl_{y_2}^{\dlt_{j_2,0}}\left[{(vu)_{[j_1,i_2]}(y_1,y_2-x_1+x_2)\over z+x_1-y_2}\right]\hspace{3.3cm}(2.39)\end{eqnarray*}
for $i_1,i_2,j_1,j_2\in\Bbb{Z}_2$ and $u\in{\cal A}_{i_1+i_2},\;v\in{\cal A}_{j_1+j_2}$ (cf. (7.3.69)-(7.3.80) in [X5]). 

Let us make a comparison of the above conformal superalgebra with the comformal algebras generating ``loop algebras" (affine Lie algebras without center) and the center-less Virasoro algebra. Let ${\cal G}$ be a Lie algebra and let $t$ be an indeterminate. Set
$$\bar{\cal G}={\cal G}\otimes_{\Bbb{F}}\Bbb{F}[t,t^{-1}]\eqno(2.40)$$
and define the algebraic operation $[\cdot,\cdot]$ on $\bar{\cal G}$ by
$$[u\otimes t^m,v\otimes t^n]=[u,v]\otimes t^{m+n}\qquad\for\;\;u,v\in {\cal G},\;m,n\in\Bbb{Z}.\eqno(2.41)$$
Then $(\bar{\cal G},[\cdot,\cdot])$ forms a Lie algebra, which is called a {\it loop Lie algebra}.\index{loop Lie algebra} 
The subspace
$${\cal B}(\bar{\cal G})={\cal G}\otimes_{\Bbb{F}}\Bbb{F}[t]\eqno(2.42)$$
forms a subalgebra of $\bar{\cal G}$. Define
$$R(\bar{\cal G})={\cal G}\otimes_{\Bbb{F}}\Bbb{F}[t^{-1}]t^{-1}\eqno(2.43)$$
and the ${\cal B}(\bar{\cal G})$-module structure on $R(\bar{\cal G})$ by
$$(u\otimes t^m)(v\otimes t^{-n})=0\;\;\mbox{if}\;\;m\geq n\;\;\mbox{and}\;\; [u,v]\otimes t^{m+n}\;\;\mbox{if}\;\;m<n\eqno(2.44)$$
for $u,v\in{\cal G},\;m\in\Bbb{N}$ and $n\in\Bbb{Z}_+$. 
Now we define a conformal algebraic structure on $R(\bar{\cal G})$ by
$$\ptl(v\otimes t^{-n})=nv\otimes t^{-n-1}\qquad\for\;\;v\in{\cal G},\;n\in\Bbb{Z}_+\eqno(2.45)$$
and
$$Y^+(u\otimes t^{-m-1},z)={1\over m!}\left({d\over dz}\right)^m(\sum_{j=0}^{\infty}(u\otimes t^j)z^{-j-1})\eqno(2.46)$$
for $u\in {\cal G}$ and $m\in\Bbb{N}$. Set
$$u[x]=\sum_{m=0}^{\infty}(u\otimes t^{-m-1})x^m\qquad\for\;\;u\in{\cal G}.\eqno(2.47)$$
Then we have
$$Y^+(u[x],z)v[y]={[u,v][y]\over z+x-y}\qquad\for\;\;u,v\in{\cal G}\eqno(2.48)$$
by (2.44) and (2.46). Moreover, $R(\bar{\cal G})$ is a free $\Bbb{F}[\ptl]$-module over ${\cal G}\otimes t^{-1}$ and the Lie algebra $\bar{\cal G}$ is generated by the conformal algebra $R(\bar{\cal G})$ (cf. Section 4.1 in [X5]). 

The {\it center-less Virasoro algebra}\index{center-less Virasoro algebra} is a vector space $\bar{\cal V}$  with a basis 
$$\{L(m)\mid m\in\Bbb{Z}\}\eqno(2.49)$$
whose Lie bracket given by
$$[L(m),L(n)]=(m-n)L(m+n)\qquad\for\;\;m,n\in\Bbb{Z}.\eqno(2.50)$$
The subspace
$${\cal B}(\bar{\cal V})=\mbox{span}\:\{L(m-1)\mid m\in\Bbb{N}\}\eqno(2.51)$$
forms a subalgebra of $\bar{\cal V}$. Define
$$R(\bar{\cal V})=\mbox{span}\:\{L(-m-2)\mid m\in\Bbb{N}\}\eqno(2.52)$$
and the ${\cal B}(\bar{\cal V})$-module structure on $R(\bar{\cal V})$ by
$$L(m-1)(L(-n-2))=\left\{\begin{array}{ll}0&\mbox{if}\;\;m\geq n+2,\\ (m+1-n)L(m-n-3)&\mbox{if}\;\;m<n+2\end{array}\right.\eqno(2.53)$$
for $m,n\in\Bbb{N}$. Now we define a conformal algebraic structure on $R(\bar{\cal V})$ by
$$\ptl(L(-n-2))=(n+1)L(-n-3)\qquad\for\;\;n\in\Bbb{N}\eqno(2.54)$$
and
$$Y^+(L(-m-2,z)={1\over m!}\left({d\over dz}\right)^m(\sum_{j=0}^{\infty}L(-j-1)z^{-j-1})\qquad\for\;\;m\in\Bbb{N}.\eqno(2.55)$$
 Set
$$L^-[x]=\sum_{m=0}^{\infty}L(-m-2)x^m.\eqno(2.56)$$
Then we have
$$Y^+(L^-[x],z)L^-[y]={dL^-(y)/dy\over z+x-y}+{2L^-(y)\over(z+x-y)^2}\eqno(2.57)$$
by (2.53) and (2.55). Moreover, $R(\bar{\cal V})$ is a free $\Bbb{F}[\ptl]$-module generated by $L(-2)$ and the Lie algebra $\bar{\cal V}$ is generated by the conformal algebra $R(\bar{V})$ (cf. Section 4.1 in [X5]).

Note that the formula in (2.39) is indeed an analogue of (2.48) and (2.57).
Let $\sgm$ be an graded involutive anti-isomorphism of ${\cal A}$, that is,
$$\sgm^2=\mbox{Id}_{\cal A},\qquad \sgm({\cal A}_i)\subset {\cal A}_i,\qquad i\in\Bbb{Z}_2,\eqno(2.58)$$
$$\sgm(e)=e,\;\;\sgm(u\cdot v)=\sgm(v)\cdot \sgm(u)\qquad\for\;\;u,v\in{\cal A}.\eqno(2.59)$$
We define
$$R({\cal A})^{\sgm}=\mbox{span}\:\{u_{[i,j]}(m,n)+(-1)^{ij}\sgm(u)_{[j,i]}(n,m)\mid u\in{\cal A},\;i,j\in\Bbb{Z}_2,\;m,n\in\Bbb{N}\}.\eqno(2.60)$$
Then $(R({\cal A})^{\sgm},\ptl,Y^+(|_{R({\cal A})^{\sgm}},z))$ forms a sub-superalgebra of $R({\cal A})$ (cf. Section 7.3 in [X5]). In next two sections, we shall prove that if ${\cal A}$ is a finite-dimensional matrix algebra, the algebras $(R({\cal A}),Y^+(\cdot,z))$ and  $(R({\cal A})^{\sgm},Y^+(\cdot,z))$ yield families of simple conformal superalgebras of finite growth.

For convenience, we shall also redenote
$$u_{[i,j]}\equiv u_{[i,j]}(0,0)\qquad\for\;\;u\in{\cal A},\;i,j\in\Bbb{Z}_2\eqno(2.61)$$
throughout Sections 3 and 4 when the context is clear. Moreover, we denote
$$Y^+(u,z)=\sum_{m=1-n}^{\infty}u(m)z^{-n-m}\qquad\for\;\;u\in R({\cal A})^{(n)},\;n\in 1+\Bbb{N}/2\eqno(2.62)$$
(cf. (2.34)).  
\psp

The following lemma will be used very often in the following two sections.
\psp

{\bf Lemma 2.1}. 
 {\it Let} $T$ {\it be a linear transformation on a vector space} $U$ {\it and let} $U_1$ {\it be a subspace of} $U$ {\it such that} $T(U_1)\subset U_1$. {\it Suppose that} $u_1,u_2,...,u_n$ {\it are eigenvectors of} $T$ {\it corresponding to different eigenvalues. If}  $\sum_{p=1}^nu_p\in U_1$, {\it then} $u_1,u_2,...,u_n\in U_1$.

\section{Simple Conformal Algebras }

In this section, we shall construct three families of infinite simple conformal algebras of finite growth from matrix algebras and prove their simplicity.

Let $k$ be a fixed positive integer. Recall that $M_{k\times k}(\Bbb{F})$ denotes the algebra of all $k\times k$-matrices with their entries in $\Bbb{F}$. Take the settings in (2.23)-(2.31) and (2.33). We let
$${\cal A}={\cal A}_0=M_{k\times k}(\Bbb{F}),\;\;{\cal A}_1=\{0\}\eqno(3.1)$$
in the general construction of $(R({\cal A}),\ptl,Y^+(\cdot,z))$. Set
$$R_{k\times k,1}=\mbox{span}\:\{u_{[1,1]}(m,n)\mid u\in{\cal A},\;m,n\in\Bbb{N}\}.\eqno(3.2)$$
and
$$R_{k\times k,\ell+2}=\mbox{span}\:\{u_{[0,0]}(m,n)\mid u\in{\cal A},\;m,n\in\Bbb{N},\;n\geq \ell\}\qquad\for\;\;\ell\in\Bbb{N}.\eqno(3.3)$$
It can be verified that all the subspaces $R_{k\times k,\ell}$ for $\ell\in\Bbb{Z}_+$ are subalgebras of $R(A)$. For each $\ell\in\Bbb{Z}_+$, the algebra $R_{k\times k,\ell}$ is $(\ell+\Bbb{N})$-weighted conformal superalgebra of finite growth with
\begin{eqnarray*} R_{k\times k,\ell}^{(n)}&=&\mbox{span}\:\{u_{[\dlt_{\ell,1},\dlt_{\ell,1}]}(m_1,m_2)\mid u\in M_{k\times k}(\Bbb{F}),\;m_1,m_2\in\Bbb{N};\\& &\qquad m_2\geq \ell-2+\dlt_{\ell,1};\;m_1+m_2+2-\dlt_{\ell,1}=n\}\hspace{4.5cm}(3.4)\end{eqnarray*}
Moreover, for convenience, we redenote
$$u(m,n)\equiv u_{[0,0]}(m,n),\;\;u\equiv u_{[0,0]}(0,0)\qquad\for\;\;u\in{\cal A}.\eqno(3.5)$$
Below, we denote by $I_k$ the $k\times k$ identity matrix.
\psp

{\bf Theorem 3.1}. {\it For each} $\ell\in\Bbb{Z}_+$, {\it the algebra} $(R_{k\times k,\ell},\ptl,Y^+(\cdot,z))$ {\it is simple. Moreover, it is generated by} $R_{k\times k,\ell}^{(\ell)}$ {\it if} $k>1$, {\it and by} $\{I_1(0,\ell),I_1(0,\ell+1)\}$ {\it (cf. (2.61) and (3.5)) if} $k=1$ {\it for } $\ell\geq 2$. {\it The algebra} $R_{k\times k,1}$ {\it is generated by} $R_{k\times k,1}^{(2)}$ {\it if} $k>1$ {\it and by} 
$$\{(I_1)_{[1,1]}(1,0),(I_1)_{[1,1]}(0,2)\}\eqno(3.6)$$
 {\it if} $k=1$.
\psp

{\it Proof}. We first consider $R_{k\times k,2+\ell}$ for $\ell\in\Bbb{N}$.  Denote by $E_{p,q}$ the matrix with 1 as its $(p,q)$-entry and $0$ as the others for $p,q\in \ol{1,k}$. Note that for $p,q,r\in \ol{1,k}$ and $m,n_1,n_2\in \Bbb{N}$ with $m,n_2\geq\ell$, we have:
$$(E_{r,r}(0,m))(0)E_{p,q}(n_1,n_2)=[\dlt_{r,p}(n_1+1)(^{-n_1-2}_{\;\;\;\;m})+\dlt_{r,q}(n_2+1)(^{n_2}_m)]E_{p,q}(n_1,n_2)\eqno(3.7)$$
by (2.30) and (2.62). Note that the coefficient in the above
$$\dlt_{r,p}(n_1+1)(^{-n_1-2}_{\;\;\;\;m})+\dlt_{r,q}(n_2+1)(^{n_2}_m)={1\over m!}{d^{m+1}\over dx^{m+1}}(\dlt_{r,q}x^{n_2+1}-\dlt_{r,p}x^{-n_1-1})|_{x=1}.\eqno(3.8)$$
If for given $p_1,p_2,q_1,q_2\in\ol{1,k}$ and $l_1,l_2,j_1,j_2\in\Bbb{N}$ with $l_2,j_2\geq \ell$,
$$\dlt_{r,p_1}(l_1+1)(^{-l_1-2}_{\;\;\;\;m})+\dlt_{r,q_1}(l_2+1)(^{l_2}_m)=\dlt_{r,p_1}(j_1+1)(^{-j_1-2}_{\;\;\;\;m})+\dlt_{r,q_1}(j_2+1)(^{j_2}_m)\eqno(3.9)$$
for any $r\in\ol{1,k}$ and $\ell\leq m\in\Bbb{N}$, then 
$$(\dlt_{r,q_1}x^{l_2+1}-\dlt_{r,p_1}x^{-l_1-1})-(\dlt_{r,q_2}x^{j_2+1}-\dlt_{r,p_1}x^{-j_1-1})\eqno(3.10)$$
is a polynomial of degree $<\ell+1$ by (3.8) and the Taylor's Theorem at $x=1$ in calculus. Since $l_2+1,j_2+1\geq \ell+1$ and $l_1,j_1\geq 0$, we have
$$(\dlt_{r,q_1}x^{l_2+1}-\dlt_{r,p_1}x^{-l_1-1})-(\dlt_{r,q_2}x^{j_2+1}-\dlt_{r,p_1}x^{-j_1-1})=0\eqno(3.11)$$
for any $r\in\ol{1,k}$. Thus
$$\mbox{(3.9) holds if and only if}\;\;p_1=p_2,\;q_1=q_2,\;j_1=l_1,\;j_2=l_2.\eqno(3.12)$$

 Let ${\cal I}$ be a nonzero ideal of $R_{k\times k,\ell+2}$ (cf. (1.5) and (1.6)). By Lemma 2.1 and (3.7)-(3.12), 
$$E_{p,q}(m,n)\in {\cal I}\qquad\mbox{for some}\;\;p,q\in\ol{1,k},\;m,n\in\Bbb{N},\;n\geq \ell.\eqno(3.13)$$
Furthermore, for $\ell\leq j\in 2\Bbb{N}$, we get 
\begin{eqnarray*}& &(E_{q,q}(0,j))(n-\ell)(E_{q,p}(0,j))(m)E_{p,q}(m,n)\\&=&(E_{q,q}(0,\ell))(n-\ell)[(m+1)(^{-m-2}_{\;\;\;\;j})E_{q,q}(0,n)+(n+1)(^{n-m}_{\;\;\;j})E_{p,p}(m,n-m)]\\&=&[(m+1)(^{-m-2}_{\;\;\;\;j})((^{-2}_{\;\:j})\dlt_{n,\ell}+(n+1)(^{\ell}_j))+\dlt_{p,q}(n+1)(^{n-m}_{\;\;\;j})((m+1)(^{-m-2}_{\;\;\;\;j})\dlt_{n,m+\ell}\\& &+(n+1)(^{\ell}_j)\dlt_{m,0})]E_{q,q}(0,\ell)\in {\cal I}\hspace{8.5cm}(3.14)\end{eqnarray*}
by (2.31) and (2.62). The coefficients of $E_{q,q}(0,\ell)$ in the last equation is positive since $j$ is even. 
Hence $E_{q,q}(0,\ell)\in {\cal I}$. 

For any $q\neq r\in\ol{1,k}$ and  $\ell< j\in 2\Bbb{N}$, we have
\begin{eqnarray*}& &(E_{r,r}(0,j))(0)(E_{r,q}(0,j))(0)(E_{q,r}(0,\ell))(0)E_{q,q}(0,\ell)\\&=&(\ell+1)(E_{r,r}(0,j))(0)(E_{r,q}(0,j))(0)
E_{q,r}(0,\ell)\\&=&(\ell+1)(j+1)(E_{r,r}(0,j))(0)E_{r,r}(0,\ell)\\&=&(\ell+1)(j+1)^2E_{r,r}(0,\ell)\in{\cal I}\hspace{9cm}(3.15)\end{eqnarray*}
by (2.31) and (2.62). Thus we have
$$I_k(0,\ell)=\sum_{r=1}^kE_{r,r}(0,\ell)\in{\cal I},\eqno(3.16)$$
where, $I_k$ is the $k\times k$ identity matrix. Moreover, (2.31) tells us that
$$(I_k(0,\ell))(-1)I_k(0,\ell)=(-1)^{\ell}(\ell+1)I_k(1,\ell)+(\ell+1)^2 I_k(0,\ell+1)\in{\cal I},\eqno(3.17)$$
$$(I_k(0,\ell+1))(-1)I_k(0,\ell)=(-1)^{\ell+1}(\ell+2)I_k(1,\ell)+(\ell+1) I_k(0,\ell+1)\in{\cal I}.\eqno(3.18)$$
Solving the above linear system, we get $I_k(0,\ell+1)\in{\cal I}$. Let $j$ be an even integer in $\{\ell,\ell+1\}$. For any $u\in {\cal A}$ (cf. (3.1)) and $m,n\in\Bbb{N}$, we have
$$(I_k(0,j))(0)u(m,\ell+n)=[(m+1)(^{-m-2}_{\;\;\;\;j})+(\ell+n+1)(^{\ell+n}_{\;\;\;j})]u(m,\ell+n)\in{\cal I}\eqno(3.19)$$
by (2.31) and (2.62). Since $j$ is even, the coefficient of $u(m,\ell+n)$ on the right-hand side is positive. Thus we get
$$u(m,\ell+n)\in{\cal I}\qquad\for\;\;u\in{\cal A},\;m,n\in\Bbb{N},\eqno(3.20)$$
that is ${\cal I}=R_{k\times k,\ell}$. So $R_{k\times k,\ell}$ is a simple conformal algebra.

Suppose that $k>1$. For $m,n\in\Bbb{N}$ and $p,q\in\ol{1,k}$ with $p\neq q$, we have
\begin{eqnarray*}& &[(E_{p,p}(0,\ell))(-1)]^m[(E_{q,q}(0,\ell))(-1)]^nE_{p,q}(0,\ell)\\&=&[\prod_{j=1}^mj(^{-j-1}_{\;\;\;\ell})\prod_{j=0}^{n-1}(\ell+j+1)(^{\ell+j}_{\;\;\ell})]E_{p,q}(m,\ell+n),\hspace{5.8cm}(3.21)\end{eqnarray*}
\begin{eqnarray*}& &(E_{p,p}(0,\ell))(0)(E_{q,p}(0,\ell))(0)E_{p,q}(m,\ell+n)\\&=&(m+1)(\ell+n+1)(^{-m-2}_{\;\;\;\;\ell})(^{\ell+n}_{\;\;\ell})E_{p,p}(m,\ell+n)\hspace{6cm}(3.22)\end{eqnarray*}
by (2.31) and (2.62). Thus $R_{k\times k,\ell+2}$ is generated by $R_{k\times k,\ell+2}^{(\ell)}$.

Let $R'$ be a subalgebra of $R_{1\times 1,\ell+2}$ generated by $\{I_1(0,\ell),I_1(0,\ell+1)\}$. We have $R_{1\times 1,\ell+2}^{(\ell+2)}=\Bbb{F}I_1(0,\ell)\subset R'$. Assume that $R_{1\times 1,\ell+2}^{(\ell+j)}\subset R'$ with $j\geq 2$. For $m,n\in\Bbb{N}$ such that $m+n=j$, we get
\begin{eqnarray*}& &(I_1(0,\ell))(-1)I_1(m,\ell+n)\\&=&(m+1)(^{-m-2}_{\;\;\;\;\;\ell})I_1(m+1,\ell+n)+(\ell+n+1)(^{\ell+n+1}_{\;\;\;\;\;\ell})I_1(m,\ell+n+1)\in R',\hspace{0.8cm}(3.23)\end{eqnarray*}
\begin{eqnarray*}&&(I_1(0,\ell+1))(-1)I_1(m,\ell+n)\\&=&(m+1)(^{-m-2}_{\;\;\:\ell+1})I_1(m+1,\ell+n)+(\ell+n+1)(^{\ell+n+1}_{\;\;\:\ell+1})I_1(m,\ell+n+1)\in R'\hspace{1cm}(3.24)\end{eqnarray*}
by (2.31) and (2.62). Since
\begin{eqnarray*}& &\left|\begin{array}{cc}(m+1)(^{-m-2}_{\;\;\;\ell}),&(\ell+n+1)(^{\ell+n+1}_{\;\;\;\;\;\ell})\\ (m+1)(^{-m-2}_{\;\;\;\ell+1}),&(\ell+n+1)(^{\ell+n+1}_{\;\;\:\ell+1})\end{array}\right|\\&=&(m+1)(\ell+n+1)(m+n+\ell+3)(\ell+1)^{-1}(^{-m-2}_{\;\;\;\;\:\ell})(^{\ell+n+1}_{\;\;\;\;\;\ell})\neq 0,\hspace{2.9cm}(3.25)\end{eqnarray*}
we have
$$I_1(m+1,\ell+n),I_1(m,\ell+n+1)\in R'.\eqno(3.26)$$
Hence $R_{1\times 1,\ell+2}^{(\ell+j+1)}\subset R'$. By induction on $j$, we have $R'=R_{k\times k,\ell+2}$.

Next we consider $R_{k\times k,1}$. Note that for $p,q,r\in \ol{1,k}$ and $m,n_1,n_2\in \Bbb{N}$, we have:
$$[(E_{r,r})_{[1,1]}(0,m)](0)(E_{p,q})_{[1,1]}(n_1,n_2)=[\dlt_{r,p}(^{-n_1-1}_{\;\;\;\;m})-\dlt_{r,q}(^{n_2}_m)](E_{p,q})_{[1,1]}(n_1,n_2)\eqno(3.27)$$
by (2.31) and (2.62). Note that the coefficient in the above
$$\dlt_{r,p}(^{-n_1-1}_{\;\;\;\;m})-\dlt_{r,q}(^{n_2}_m)={1\over m!}{d^m\over dx^m}(\dlt_{r,p}x^{-n_1-1}-\dlt_{r,q}x^{n_2})|_{x=1}.\eqno(3.28)$$

 Let ${\cal I}$ be a nonzero ideal of $R_{k\times k,1}$ (cf. (1.5) and (1.6)). By Lemma 2.1, (3.27), (3.28) and the Taylor's theorem at $x=1$ in calculus (cf. (3.9)-(3.12)), 
$$(E_{p,q})_{[1,1]}(m,n)\in {\cal I}\qquad\mbox{for some}\;\;p,q\in\ol{1,k},\;m,n\in\Bbb{N}.\eqno(3.29)$$
Let
$$R^{\dg}=\mbox{span}\{u_{[1,1]}(j,l)\mid u\in{\cal A},\;j,l\in\Bbb{N},\;l\geq 1\},\eqno(3.30)$$
$$R^{\ast}=\mbox{span}\{u_{[1,1]}(l,j)\mid u\in{\cal A},\;j,l\in\Bbb{N},\;l\geq 1\}.\eqno(3.31)$$
Then $(R^{\dg},\ptl,Y^+(\cdot,z))$ and $(R^{\ast},\ptl,Y^+(\cdot,z))$ are both subalgebras of $R_{k\times k,1}$, which are isomorphic to $(R_{k\times k,2},\ptl,Y^+(\cdot,z))$ through the following correspondences:
$$-(l+1)u_{[1,1]}(j,l+1)\leftrightarrow u_{[0,0]}(j,l),\;\;(l+1)u_{[1,1]}(l+1,j)\leftrightarrow u_{[0,0]}(l,j)\eqno(3.32)$$
for $u\in{\cal A}$ and $j,l\in\Bbb{N}$. The correspondences in the above are {\it boson-fermion correspondences}\index{boson-fermion correspondence} in physics. 
If $m+n>0$ in (3.29), we have
$$(E_{p,q})_{[1,1]}(m,n)\in R^{\dg}\bigcup R^{\ast},\eqno(3.33)$$
which implies 
$$R^{\dg}\subset{\cal I}\;\;\mbox{or}\;\;R^{\ast}\subset {\cal I}.\eqno(3.34)$$
In particular, we have
$$(I_k)_{[1,1]}(0,1)\in {\cal I}\;\;\mbox{or}\;\;(I_k)_{[1,1]}(1,0)\in{\cal I}.\eqno(3.35)$$
Assume that $m=n=0$ in (3.29). Then we have
$$[(I_k)_{[1,1]}(0,1)](-1)(E_{p,q})_{[1,1]}=-(E_{p,q})_{[1,1]}(1,0)-(E_{p,q})_{[1,1]}(0,1)\in{\cal I}\eqno(3.36)$$
(cf. (2.61)). Thus we can always assume $m+n>0$ by Lemma 2.1, (3.27) and (3.28). So (3.33) holds. Furthermore,
$$[(I_k)_{[1,1]}(0,1)](0)u_{[1,1]}(j,l)=-(j+l+1)u_{[1,1]}(j,l)\qquad\for\;\;u\in{\cal A},\;j,l\in\Bbb{N},\eqno(3.37)$$
$$[(I_k)_{[1,1]}(1,0)](0)u_{[1,1]}(j,l)=(j+l+1)u_{[1,1]}(j,l)\qquad\for\;\;u\in{\cal A},\;j,l\in\Bbb{N}\eqno(3.38)$$
by (2.31) and (2.62). Therefore, ${\cal I}=R_{k\times k,1}$, that is, $R_{k\times k,1}$ is simple.

Assume that $k>1$. Let $R'$ be the subalgebra of $R_{k\times k,1}$ generated by $R_{k\times k,1}^{(2)}$.  Then
$R^{\dg},R^{\ast}\in R'$ by the fact that $R_{k\times k,2}$ is generated by $R_{k\times k,2}^{(2)}$ and the isomorphisms in 
(3.32). Note that
$$[(I_k)_{[1,1]}(1,0)](1)(E_{p,q})_{[1,1]}(0,1)=2(E_{p,q})_{[1,1]}(0,0)\in R'\qquad\for\;\;p,q\in\ol{1,k}\eqno(3.39)$$
by (2.31) and (2.62). So $R_{k\times k,1}=R'$. 

Assume that $k=1$. Let $R'$ be the subalgebra of $R_{1\times 1,1}$ generated by (3.6). Note that
$$[(I_1)_{[1,1]}(1,0)](1)(I_1)_{[1,1]}(0,2)=3(I_1)_{[1,1]}(0,1)\in R',\eqno(3.40)$$
$$[(I_1)_{[1,1]}(0,2)](-1)(I_1)_{[1,1]}(1,0)=6(I_1)_{[1,1]}(2,0)\eqno(3.41)$$
by (2.31) and (2.62). By (3.39), the fact $R_{1\times 1,2}$ is generated by $\{I_1,I_1(0,1)\}$ and the isomorphisms in 
(3.32), we have $R'=R_{1\times 1,1}.\qquad\Box$
\psp

Let $\sgm_1:A\mapsto A^t$ be the transpose map of matrices. Then $\sgm_1$ is an involutive anti-isomorphism of $M_{k\times k}(\Bbb{F})$. Thus we have the following subalgebras of $R({\cal A})^{\sgm}$ (cf. (2.67)):
$$R^{\ast}_{k\times k,1}=\mbox{span}\:\{(E_{p,q})_{[1,1]})(m,n)-(E_{q,p})_{[1,1]})(n,m)\mid  p,q\in\ol{1,k},\;m,n\in \Bbb{N}\},\eqno(3.42)$$
$$R^{\ast}_{k\times k,2}=\mbox{span}\:\{E_{p,q}(m,n)+E_{q,p}(n,m)\mid p,q\in\ol{1,k},\;m,n\in \Bbb{N}\}.\eqno(3.43)$$
\psp

{\bf Theorem 3.2}.  {\it The algebras} $(R^{\ast}_{k\times k,1},\ptl,Y^+(\cdot,z))$ {\it and} $(R^{\ast}_{k\times k,2},\ptl,Y^+(\cdot,z))$ {\it are simple. Moreover,} $R^{\ast}_{k\times k,1}$ {\it is generated by} 
$$(R^{\ast}_{k\times k,1})^{(2)}=\mbox{span}\:\{(E_{p,q})_{[1,1]})(1,0)-(E_{q,p})_{[1,1]})(0,1)\mid  p,q\in\ol{1,k},\;m,n\in \Bbb{N}\}\eqno(3.44)$$
{\it when} $k>1$ {\it and by} 
$$\{(I_1)_{[1,1]}(1,0)-(I_1)_{[1,1]}(0,1),(I_1)_{[1,1]}(3,0)-(I_1)_{[1,1]}(0,3)\}\eqno(3.45)$$
{\it if} $k=1$. {\it The algebra} $R^{\ast}_{k\times k,2}$
{\it is generated by} 
$$(R^{\ast}_{k\times k,2})^{(2)}=\mbox{span}\:\{E_{p,q}+E_{q,p}\mid  p,q\in\ol{1,k},\;m,n\in \Bbb{N}\}\eqno(3.46)$$
{\it (cf. (3.5)) when} $k>1$ {\it and by}
$\{I_1,I_1(1,1)\}$ {\it if} $k=1$. 
\psp

{\it Proof}. For $p,q,r\in\ol{1,k}$ and $m,n_1,n_2\in\Bbb{N}$, we have
\begin{eqnarray*}& &[(E_{r,r})_{[1,1]}(m,0)-(E_{r,r})_{[1,1]}(0,m)](0)[(E_{p,q})_{[1,1]}(n_1,n_2)-(E_{q,p})_{[1,1]}(n_2,n_1)]\\&=&[\dlt_{r,p}[(^{n_1}_m)-(^{-n_1-1}_{\;\;\;\;m})]+\dlt_{r,q}[(^{n_2}_m)-(^{-n_2-1}_{\;\;\;\;m})]]\\& &[(E_{p,q})_{[1,1]}(n_1,n_2)-(E_{q,p})_{[1,1]}(n_2,n_1)],\hspace{7.2cm}(3.47)\end{eqnarray*}
\begin{eqnarray*}& &(E_{r,r})(m,0)+E_{r,r}(0,m))(0)(E_{p,q}(n_1,n_2)+E_{q,p}(n_2,n_1))\\&=&[\dlt_{r,p}(n_1+1)[(^{n_1}_m)+(^{-n_1-2}_{\;\;\;\;m})]+\dlt_{r,q}(n_2+1)[(^{n_2}_m)+(^{-n_2-2}_{\;\;\;\;m})]]\\& &[E_{p,q}(n_1,n_2)+E_{q,p}(n_2,n_1)]\hspace{9.2cm}(3.48)\end{eqnarray*}
by (2.31) and (2.62). Note that
\begin{eqnarray*}& &\dlt_{r,p}[(^{n_1}_m)-(^{-n_1-1}_{\;\;\;\;m})]+\dlt_{r,q}[(^{n_2}_m)-(^{-n_2-1}_{\;\;\;\;m})]\\&=&{1\over m!}{d^m\over dx^m}[\dlt_{r,p}x^{n_1}+\dlt_{r,q}x^{n_2}-\dlt_{r,p}x^{-n_1-1}-\dlt_{r,q}x^{-n_2-1}]|_{x=1},\hspace{4.2cm}(3.49)\end{eqnarray*}
\begin{eqnarray*}& &\dlt_{r,p}(n_1+1)[(^{n_1}_m)+(^{-n_1-2}_{\;\;\;\;m})]+\dlt_{r,q}(n_2+1)[(^{n_2}_m)+(^{-n_2-2}_{\;\;\;\;m})]\\&=&{1\over m!}{d^{m+1}\over dx^{m+1}}[\dlt_{r,p}x^{n_1+1}+\dlt_{r,q}x^{n_2+1}-\dlt_{r,p}x^{-n_1-2}-\dlt_{r,q}x^{-n_2-2}]|_{x=1}.\hspace{3.1cm}(3.50)\end{eqnarray*}

Let ${\cal I}$ be a nonzero ideal of $R^{\ast}_{k\times k,1}$. By Lemma 2.1, (3.47), (3.49) and the Taylor's theorem at $x=1$ in calculus (cf. (3.9)-(3.12)), we have
$$(E_{p,q})_{[1,1]}(n_1,n_2)-(E_{q,p})_{[1,1]}(n_2,n_1)\in{\cal I}\qquad\mbox{for some}\;\;p,q\in\ol{1,k},\;n_1,n_2\in\Bbb{N}\eqno(3.51)$$
such that $p\neq q$ or $n_1\neq n_2$. If $p\neq q$, we have
\begin{eqnarray*}& &[(E_{q,q})_{[1,1]}(1,0)-(E_{q,q})_{[1,1]}(0,1)](n_2)[(E_{p,q})_{[1,1]}(1,0)-(E_{q,p})_{[1,1]}(0,1)](n_1-1)\\& &[(E_{p,q})_{[1,1]}(n_1,n_2)-(E_{q,p})_{[1,1]}(n_2,n_1)]\\&=&[(E_{q,q})_{[1,1]}(1,0)-(E_{q,q})_{[1,1]}(0,1)](n_2)[(n_1-n_2-1)[(E_{p,p})_{[1,1]}(n_1,n_2+1-n_1)\\& &-(E_{p,p})_{[1,1]}(n_2+1-n_1,n_1)]+(n_1+1)[(E_{q,q})_{[1,1]}(1,n_2)-(E_{q,q})_{[1,1]}(n_2,1)]]\\&=&(n_1+1)(3\dlt_{n_2,0}+n_2+1-2\dlt_{n_2,1})[(E_{q,q})_{[1,1]}(1,0)-(E_{q,q})_{[1,1]}(0,1)],\hspace{2.2cm}(3.52)\end{eqnarray*}
\begin{eqnarray*}& &[(E_{q,q})_{[1,1]}(1,0)-(E_{q,q})_{[1,1]}(0,1)](n_2-1)[(E_{p,q})_{[1,1]}(1,0)-(E_{q,p})_{[1,1]}(0,1)](n_1)\\& &[(E_{p,q})_{[1,1]}(n_1,n_2)-(E_{q,p})_{[1,1]}(n_2,n_1)]\\&=&[(E_{q,q})_{[1,1]}(1,0)-(E_{q,q})_{[1,1]}(0,1)](n_2-1)[(n_2-n_1)[(E_{p,p})_{[1,1]}(n_1,n_2-n_1)\\& &-(E_{p,p})_{[1,1]}(n_2-n_1,n_1)]+(n_1+1)[(E_{q,q})_{[1,1]}(0,n_2)-(E_{q,q})_{[1,1]}(n_2,0)]]\\&=&(n_1+1)(2\dlt_{0,n_2}-n_2-2-\dlt_{1,n_2})((E_{q,q})_{[1,1]}(1,0)-(E_{q,q})_{[1,1]}(0,1))\in{\cal I}\hspace{1.6cm}(3.53)\end{eqnarray*}
by (2.31) and (2.62). Since (3.52) is zero only if $n_2=1$ and (3.53) is zero only if $n_2=0$, we have 
$$(E_{q,q})_{[1,1]}(1,0)-(E_{q,q})_{[1,1]}(0,1)\in{\cal I}.\eqno(3.54)$$
Assume $p=q$ and $n_1\neq n_2$. We have
\begin{eqnarray*}& &[(E_{q,q})_{[1,1]}(1,0)-(E_{q,q})_{[1,1]}(0,1)](n_2-1)[(E_{q,q})_{[1,1]}(1,0)-(E_{q,q})_{[1,1]}(0,1)](n_1)\\& &[(E_{q,q})_{[1,1]}(n_1,n_2)-(E_{q,q})_{[1,1]}(n_2,n_1)]\\ &=& [(E_{q,q})_{[1,1]}(1,0)-(E_{q,q})_{[1,1]}(0,1)](n_2-1)[(2n_2-n_1+1)[(E_{q,q})_{[1,1]}(n_1,n_2-n_1)\\& &-(E_{q,q})_{[1,1]}(n_2-n_1,n_1)]+(n_1+1)[(E_{q,q})_{[1,1]}(0,n_2)-(E_{q,q})_{[1,1]}(n_2,0)]]\\&=&(2n_2-n_1+1)[(n_2-2n_1+2)[(E_{q,q})_{[1,1]}(n_1,1-n_1)-
(E_{q,q})_{[1,1]}(1-n_1,n_1)]\\& &+(n_1+1)\dlt_{n_1,n_2+1}[(E_{q,q})_{[1,1]}(0,1)-(E_{q,q})_{[1,1]}(1,0)]]\\&&
+(n_1+1)\dlt_{n_1+1,n_2}(n_2+2-2\dlt_{n_2,0}+\dlt_{n_2,1})[(E_{q,q})_{[1,1]}(0,1)-(E_{q,q})_{[1,1]}(1,0)]\\&=&\mu[(E_{q,q})_{[1,1]}(0,1)-(E_{q,q})_{[1,1]}(1,0)],\hspace{7.7cm}(3.55)\end{eqnarray*}
for some $0\neq\mu\in\Bbb{F}$, and $\mu=0$ only if $(n_1,n_2)\in \{(1,0),(1,2),(2+\Bbb{N},0)\}$. So (3.54) holds if $(n_1,n_2)\not\in\{(1,0),(1,2),(2+\Bbb{N},0)\}$. Symmetrically, we can prove (3.54) when $(n_2,n_1)\not\in\{(1,0),(1,2),(2+\Bbb{N},0)\}$. Thus (3.54) always holds.

Let $q\neq j\in\ol{1,k}$. Observe that
\begin{eqnarray*}& &[[(E_{j,q})_{[1,1]}(0,1)-(E_{q,j})_{[1,1]}(1,0)](0)]^2[(E_{q,q})_{[1,1]}(1,0)-(E_{q,q})_{[1,1]}(0,1)]\\&=&[(E_{j,q})_{[1,1]}(0,1)-(E_{q,j})_{[1,1]}(1,0)](0)[(E_{j,q})_{[1,1]}(0,1)-(E_{q,j})_{[1,1]}(1,0)\\& &+2((E_{q,j})_{[1,1]}(0,1)-(E_{j,q})_{[1,1]}(1,0))]\\&=&4[(E_{j,j})_{[1,1]}(1,0)-(E_{j,j})_{[1,1]}(0,1)]+2[(E_{q,q})_{[1,1]}(1,0)-(E_{q,q})_{[1,1]}(0,1)]\in{\cal I}\hspace{0.8cm}(3.56)\end{eqnarray*}
by (2.31) and (2.62). Thus $(E_{j,j})_{[1,1]}(1,0)-(E_{j,j})_{[1,1]}(0,1)\in {\cal I}$, and  (3.54) holds for any $q\in\ol{1,k}$. Hence
$$(I_k)_{[1,1]}(1,0)-(I_k)_{[1,1]}(0,1)=\sum_{q=1}^k[(E_{q,q})_{[1,1]}(1,0)-(E_{q,q})_{[1,1]}(0,1)]\in{\cal I}.\eqno(3.57)$$
By (3.37) and (3.38), ${\cal I}=R_{k\times k,1}$. So $R_{k\times k,1}$ is simple.

Assume  $k>1$. For $m,n\in\Bbb{N}$ and $p,q\in\ol{1,k}$ such that $p\neq q$, we have
\begin{eqnarray*}& &([(E_{p,p})_{[1,1]}(1,0)-(E_{p,p})_{[1,1]}(0,1)](-1))^m([(E_{q,q})_{[1,1]}(1,0)-(E_{q,q})_{[1,1]}(0,1)](-1))^n\\& &[(E_{p,q})_{[1,1]}(1,0)-(E_{q,p})_{[1,1]}(0,1)]\\&=&2^{m+n}(m+1)!n![(E_{p,q})_{[1,1]}(m+1,n)-(E_{q,p})_{[1,1]}(n,m+1)],\hspace{3.4cm}(3.58)\end{eqnarray*}
\begin{eqnarray*}& &[(E_{p,p})_{[1,1]}(1,0)-(E_{p,p})_{[1,1]}(0,1)](1)[(E_{p,q})_{[1,1]}(1,0)-(E_{q,p})_{[1,1]}(0,1)]\\&=&2[(E_{p,q})_{[1,1]}
-(E_{q,p})_{[1,1]}]\hspace{9.8cm}(3.59)\end{eqnarray*}
by (2.31) and (2.62). Moreover, we have
\begin{eqnarray*}&&[(E_{p,p})_{[1,1]}(1,0)-(E_{p,p})_{[1,1]}(0,1)](0)[(E_{p,q})_{[1,1]}(0,1)-(E_{q,p})_{[1,1]}(1,0)](0)\\& &[(E_{p,q})_{[1,1]}(m,n)-(E_{q,p})_{[1,1]}(n,m)]\\&=&[(E_{p,p})_{[1,1]}(1,0)-(E_{p,p})_{[1,1]}(0,1)](0)\{(n+1)[(E_{p,p})_{[1,1]}(n,m)-(E_{p,p})_{[1,1]}(m,n)]\\& &+m[(E_{q,q})_{[1,1]}(n,m)-(E_{q,q})_{[1,1]}(m,n)]\}\\&=&2(n+1)(m+n+1)[(E_{p,p})_{[1,1]}(n,m)-(E_{p,p})_{[1,1]}(m,n)]\hspace{4.1cm}(3.60)\end{eqnarray*}
for $m,n\in\Bbb{N}$ and $p,q\in\ol{1,k}$ with $p\neq q$. Thus $R_{k\times k,1}^{\ast}$ is generated $(R_{k\times k,1}^{\ast})^{(2)}$.

Now we assume $k=1$. Let $R'$ be the subalgebra of $R_{k\times k,1}^{\ast}$ generated by (3.45). By (3.42), $(R_{k\times k,1}^{\ast})^{(1)}=\{0\}$ and $(R_{k\times k,1}^{\ast})^{(2)}=\Bbb{F}((I_1)_{[1,1]}(1,0)-(I_1)_{[1,1]}(0,1))\subset R'$. Assume that $(R_{k\times k,1}^{\ast})^{(j)}\subset R'$ for some $2\leq j\in\Bbb{N}$. For $m,n\in\Bbb{N}$ such that $m+n=j$, we get
\begin{eqnarray*}& &[(I_1)_{[1,1]}(1,0)-(I_1)_{[1,1]}(0,1)](-1)[(I_1)_{[1,1]}(m,n)-(I_1)_{[1,1]}(n,m)]\\&=&2(m+1)[(I_1)_{[1,1]}(m+1,n)-(I_1)_{[1,1]}(n,m+1)]\\&&+2(n+1)[(I_1)_{[1,1]}(m,n+1)-(I_1)_{[1,1]}(n+1,m)],\hspace{5.1cm}(3.61)\end{eqnarray*}
\begin{eqnarray*}& &[(I_1)_{[1,1]}(3,0)-(I_1)_{[1,1]}(0,3)](-1)[(I_1)_{[1,1]}(m,n)-(I_1)_{[1,1]}(n,m)]\\&=&[(^{m+1}_{\;\;\;3})-(^{-m-1}_{\;\;\;\;3})][(I_1)_{[1,1]}(m+1,n)-(I_1)_{[1,1]}(n,m+1)]\\& &+[(^{n+1}_{\;\;3})-(^{-n-1}_{\;\;\;\;3})][(I_1)_{[1,1]}(m,n+1)-(I_1)_{[1,1]}(n+1,m)]\\&=&{1\over 3}\{(m+1)(m^2+2m+3)[(I_1)_{[1,1]}(m+1,n)-(I_1)_{[1,1]}(n,m+1)]\\& &+(n+1)(n^2+2n+3)[(I_1)_{[1,1]}(m,n+1)-(I_1)_{[1,1]}(n+1,m)]\}\hspace{2.8cm}(3.62)\end{eqnarray*}
by (2.31) and (2.62). When $m\neq n$, we have
\begin{eqnarray*}& &\left|\begin{array}{cc}m+1,&n+1\\ (m+1)(m^2+2m+3),&(n+1)(n^2+2n+3)\end{array}\right|\\&=&(n-m)(m+n+2)(m+1)(n+1)\neq 0.\hspace{6.7cm}(3.63)\end{eqnarray*}
Thus 
$$(I_1)_{[1,1]}(m+1,n)-(I_1)_{[1,1]}(n,m+1),(I_1)_{[1,1]}(m,n+1)-(I_1)_{[1,1]}(n+1,m)\in R'\eqno(3.64)$$
if $m\neq n$. Since
$$(n,n+1)=(n-1,n+1)+(1,0),\;\;(n+1,n)=(n+1,n-1)+(0,1),\eqno(3.65)$$
we have $(R_{k\times k,1}^{\ast})^{(j+1)}\subset R'$. By induction on $j$, $R_{k\times k,1}^{\ast}=R'$.

Next we consider $R_{k\times k,2}^{\ast}$. Let ${\cal I}$ be a nonzero ideal of $R_{k\times k,2}^{\ast}$. By Lemma 2.1, (3.48), (3.50) and the Taylor's theorem at $x=1$ in calculus (cf. (3.9)-(3.12)), we have
$$E_{p,q}(m,n)+E_{q,p}(n,m)\in{\cal I}\qquad\mbox{for some}\;\;p,q\in\ol{1,k},\;m,n\in\Bbb{N}.\eqno(3.66)$$
Moreover, we get
\begin{eqnarray*}& &E_{q,q}(n)(E_{p,q}+E_{q,p})(m)(E_{p,q}(m,n)+E_{q,p}(n,m))\\&=&E_{q,q}(n)[(m+1)(E_{q,q}(0,n)+E_{q,q}(n,0))\\& &+(n+1)(E_{p,p}(m,n-m)+E_{p,p}(n-m,m))]\\&=&2(m+1)(n+1)(1+\dlt_{0,n})E_{q,q}\in{\cal I}\hspace{7.9cm}(3.67)\end{eqnarray*}
if $p\neq q$ and
\begin{eqnarray*}& &E_{q,q}(n)E_{q,q}(m)(E_{q,q}(m,n)+E_{q,q}(n,m))\\&=&E_{q,q}(n)[(m+1)(E_{q,q}(0,n)+E_{q,q}(n,0))\\& &+(n+1)(E_{q,q}(n-m,m)+E_{q,q}(m,n-m))]\\&=&2(m+1)(n+1)[(1+\dlt_{n,0})(1+(\dlt_{m,n}+\dlt_{m,0})/(1+\dlt_{m,0}\dlt_{n,0})]E_{q,q}\in{\cal I}\hspace{2cm}(3.68)\end{eqnarray*}
by (2.31), (2.62) and (3.5). So $E_{q,q}\in {\cal I}$. Let $q\neq j\in\ol{1,k}$. We have
$$(E_{j,q}+E_{q,j})_1^2(E_{q,q})=2(E_{j,j}+E_{q,q})\in {\cal I}\eqno(3.69)$$
by (2.31) and (2.62). This implies $E_{j,j}\in{\cal I}$. Thus $I_k=\sum_{i=1}^kE_{i,i}\in {\cal I}$. By (3.19) with $j=\ell=0$, ${\cal I}=R^{\ast}_{k\times k,2}$. Hence $R_{k\times k,2}^{\ast}$ is simple.

 Assume $k>1$. For any $m,n\in\Bbb{N}$ and $p,q\in\ol{1,k}$ such that $p\neq q$, we have
\begin{eqnarray*}& &(E_{p,p}(-1))^m(E_{q,q}(-1))^n(E_{p,q}+E_{q,p})\\&=&(m+1)!(n+1)!(E_{p,q}(m,n)+E_{q,p}(n,m),\hspace{6.6cm}(3.70)\end{eqnarray*}
\begin{eqnarray*}& &E_{p,p}(0)(E_{p,q}+E_{q,p})(0)(E_{p,q}(m,n)+E_{q,p}(n,m))\\&=&(m+n+2)(n+1)(E_{p,p}(m,n)+E_{p,p}(n,m))\hspace{6.1cm}(3.71)\end{eqnarray*}
by (2.31) and (2.62). So $R_{k\times k,2}^{\ast}$ is generated by $(R_{k\times k,2}^{\ast})^{(2)}$.

Now we assume $k=1$. Let $R'$ be the subalgebra of $R_{k\times k,2}^{\ast}$ generated by $\{I_1,I_1(1,1)\}$. Note that
\begin{eqnarray*}& &I_1(-1)(I_1(m,n)+I_1(n,m))=(m+1)(I_1(m+1,n)\\& &+I_1(n,m+1))+(n+1)(I_1(m,n+1)+I_1(n+1,m)),\hspace{4.5cm}(3.72)\end{eqnarray*}
\begin{eqnarray*}\hspace{1cm}& &(I_1(1,1))(-1)(I_1(m,n)+I_1(n,m))\\&=&-(m+1)^2(m+2)(I_1(m+1,n)+I_1(n,m+1))\\&&-(n+1)^2(n+2)(I_1(m,n+1)+I_1(n+1,m)\hspace{4.9cm}(3.73)\end{eqnarray*}
for $m,n\in\Bbb{N}$. Since
$$\left|\begin{array}{cc}m+1,&n+1\\-(m+1)^2(m+2),&-(n+1)^2(n+2)\end{array}\right|=(m-n)(m+n+3)(m+1)(n+1),\eqno(3.74)$$
we can prove that $I_1(m,n)\in R'$ by mathematical induction on $m+n$. Therefore, we have
$R^{\ast}_{k\times k,2}=R'. \qquad\Box$
\psp

Assume that $k=2k_1$ is an even integer. For
$$A=\left(\begin{array}{cc}A_{1,1},& A_{1,2}\\ A_{2,1},& A_{2,2}\end{array}\right)\;\;\mbox{with}\;\;A_{i,j}\in M_{k_1\times k_1}(\Bbb{F}),\eqno(3.75)$$
we define
$$\sgm_2(A)=\left(\begin{array}{cc}&I_{k_1}\\-I_{k_1}&\end{array}\right)A^t\left(\begin{array}{cc}&-I_{k_1}\\ I_{k_1}&\end{array}\right)=\left(\begin{array}{cc}A_{2,2}^t,&-A_{1,2}^t\\ -A_{2,1}^t,&A_{1,1}^t\end{array}\right),\eqno(3.76)$$
where the empty entries are zero. Then $\sgm_2$ is another involutive anti-isomorphism of $M_{k\times k}(\Bbb{F})$. Moreover, we have the following subalgebra of $R_{k\times k,1}$:
\begin{eqnarray*}R^{\dg}_{k\times k,1}&=&\mbox{span}\;\{(E_{p,q})_{[1,1]}(m,n)-(E_{k_1+q,+k_1+p})_{[1,1]}(n,m),\\& &(E_{p,k_1+q})_{[1,1]}(m,n)+(E_{q,k_1+p})_{[1,1]}(n,m)\\&& (E_{k_1+p,q})_{[1,1]}(m,n)+(E_{q,k_1+p})_{[1,1]}(n,m)\mid p,q\in\ol{1,k_1},\;m,n\in\Bbb{N}\}\hspace{1.4cm}(3.77)\end{eqnarray*}
and the subalgebra of $R_{k\times k,2}$:
\begin{eqnarray*}R^{\dg}_{k\times k,2}&=&\mbox{span}\;\{E_{p,q}(m,n)+E_{k_1+q,+k_1+p}(n,m),\;E_{p,k_1+q}(m,n)-E_{q,k_1+p}(n,m),\\&& E_{k_1+p,q}(m,n)-E_{q,k_1+p}(n,m)\mid p,q\in\ol{1,k_1},\;m,n\in\Bbb{N}\}.\hspace{3.2cm}(3.78)\end{eqnarray*}

{\bf Theorem 3.3}.  {\it The conformal algebras} $(R^{\dg}_{k\times k,1},\ptl,Y^+(\cdot,z))$ {\it and} $(R^{\dg}_{k\times k,2},\ptl,Y^+(\cdot,z))$ {\it are simple. Moreover, the algebra} $(R^{\dg}_{k\times k,1},\ptl,Y^+(\cdot,z))$ {\it is generated by}
\begin{eqnarray*}(R^{\dg}_{k\times k,1})^{(2)}&=&\mbox{span}\;\{(E_{p,q})_{[1,1]}(\es,1-\es)-(E_{k_1+q,+k_1+p})_{[1,1]}(1-\es,\es),\\&&(E_{p,k_1+q})_{[1,1]}(1,0)+(E_{q,k_1+p})_{[1,1]}(0,1),(E_{k_1+p,q})_{[1,1]}(1,0)\\& & +(E_{q,k_1+p})_{[1,1]}(0,1)\mid p,q\in\ol{1,k_1},\;\es=0,1\}\hspace{4.6cm}(3.79)\end{eqnarray*}
{\it when} $k_1>1$ {\it and by}
$$\begin{array}{l}\{(E_{1,2})_{[1,1]}(2,0)+(E_{1,2})_{[1,1]}(0,2),\;(E_{2,1})_{[1,1]}(2,0)+(E_{2,1})_{[1,1]}(0,2),\\(E_{1,1})_{[1,1]}(0,1)+(E_{1,2})_{[1,1]}(1,0)\}\end{array}\eqno(3.80)$$
{\it if} $k_1=1$. {\it The algebra} $(R^{\dg}_{k\times k,2},\ptl,Y^+(\cdot,z))$
 {\it is generated by} 
\begin{eqnarray*}(R_{k\times k,2}^{\dg})^{(2)}&=&\mbox{span}\:\{E_{p,q}+E_{k_1+q,+k_1+p},\;E_{p,k_1+q}-E_{q,k_1+p},\\& &E_{k_1+p,q}-E_{q,k_1+p}\mid p,q\in\ol{1,k_1}\}\hspace{6.6cm}(3.81)\end{eqnarray*}
{\it if} $k_1>1$ {\it and by}
 $$\{I_2,\;E_{1,2}(1,0)-E_{1,2}(0,1),\;E_{2,1}(1,0)-E_{2,1}(0,1)\}\eqno(3.82)$$
 {\it if} $k=1$ {\it (cf. (3.5))}.
\psp

{\it Proof}. For $p,q,r\in\ol{1,k_1}$ and $m,n_1,n_2\in\Bbb{N}$, we have
\begin{eqnarray*}& &[(E_{r,r})_{[1,1]}(0,m)-(E_{k_1+r,k_1+r})_{[1,1]}(m,0)](0)\\& &[(E_{p,q})_{[1,1]}(n_1,n_2)-(E_{k_1+q,k_1+p})_{[1,1]}(n_2,n_1)]\\&=&[\dlt_{r,p}(^{-n_1-1}_{\;\;\;\;m})-\dlt_{r,q}(^{n_2}_m)][(E_{p,q})_{[1,1]}(n_1,n_2)-(E_{k_1+q,k_1+p})_{[1,1]}(n_2,n_1)],\hspace{2.3cm}(3.83)\end{eqnarray*}
\begin{eqnarray*}& &[(E_{r,r})_{[1,1]}(0,m)-(E_{k_1+r,k_1+r})_{[1,1]}(m,0)](0)\\& &[(E_{p,k_1+q})_{[1,1]}(n_1,n_2)+(E_{q,k_1+p})_{[1,1]}(n_2,n_1)]\\&=&[\dlt_{r,p}(^{-n_1-1}_{\;\;\;\;m})+\dlt_{r,q}(^{-n_2-1}_{\;\;\;\;m})][(E_{p,k_1+q})_{[1,1]}(n_1,n_2)+(E_{q,k_1+p})_{[1,1]}(n_2,n_1)],\hspace{1.7cm}(3.84)\end{eqnarray*}
\begin{eqnarray*}& &[(E_{r,r})_{[1,1]}(0,m)-(E_{k_1+r,k_1+r})_{[1,1]}(m,0)](0)\\& &[(E_{k_1+p,q})_{[1,1]}(n_1,n_2)+(E_{k_1+q,p})_{[1,1]}(n_2,n_1)]\\&=&[-\dlt_{r,p}(^{n_1}_m)-\dlt_{r,q}(^{n_2}_m)][(E_{k_1+p,q})_{[1,1]}(n_1,n_2)+(E_{k_1+q,p})_{[1,1]}(n_2,n_1)]\hspace{2.8cm}(3.85)\end{eqnarray*}
by (2.31) and (2.62). Moreover,
$$\dlt_{r,p}(^{-n_1-1}_{\;\;\;\;m})-\dlt_{r,q}(^{n_2}_m)={1\over m!}{d^m\over dx^m}(\dlt_{r,p}x^{-n_1-1}-\dlt_{r,q}x^{n_2})|_{x=1},\eqno(3.86)$$
$$\dlt_{r,p}(^{-n_1-1}_{\;\;\;\;m})+\dlt_{r,q}(^{-n_2-1}_{\;\;\;\;m})={1\over m!}{d^m\over dx^m}(\dlt_{r,p}x^{-n_1-1}+\dlt_{r,q}x^{-n_2-1})|_{x=1},\eqno(3.87)$$
$$-\dlt_{r,p}(^{n_1}_m)-\dlt_{r,q}(^{n_2}_m)={1\over m!}{d^m\over dx^m}(-\dlt_{r,p}x^{n_1}-\dlt_{r,q}x^{n_2})|_{x=1}.\eqno(3.88)$$

Let ${\cal I}$ be a nonzero ideal of $R^{\dg}_{k\times k,1}$. By Lemma 2.1, (3.83)-(3.88) and the Taylor's theorem at $x=1$ in calculus (cf. (3.9)-(3.12)),  ${\cal I}$ contains at least one of the following elements:
$$\begin{array}{l}\{(E_{p,q})_{[1,1]}(n_1,n_2)-(E_{k_1+q,k_1+p})_{[1,1]}(n_2,n_1),\\ (E_{p,k_1+q})_{[1,1]}(n_1,n_2)+(E_{q,k_1+p})_{[1,1]}(n_2,n_1),\\
(E_{k_1+p,q})_{[1,1]}(n_1,n_2)+(E_{k_1+q,p})_{[1,1]}(n_2,n_1)\}\end{array}\eqno(3.89)$$
for some $p,q\in\ol{1,k_1}$ and $n_1,n_2\in\Bbb{N}$. Note that the subspace
$$\bar{R}=\mbox{span}\{(E_{j_1,j_2})_{[1,1]}(m,n)-(E_{k_1+j_2,k_1+j_1})_{[1,1]}(n,m)\mid j_1,j_2\in\ol{1,k_1},\;m,n\in\Bbb{N}\}\eqno(3.90)$$
forms a subalgebra of $R^{\dg}_{k\times k,1}$ that is isomorphic to $R_{k_1\times k_1,1}$. By Theorem 3.1, $R_{k_1\times k_1,1}$ is simple. Hence
$$\bar{R}\in{\cal I}\eqno(3.91)$$
if the first element in (3.89) is in ${\cal I}$. Assume that the second element in (3.89) is in ${\cal I}$. Using (2.61), we have
\begin{eqnarray*}& &[(E_{k_1+p,q})_{[1,1]}+(E_{k_1+q,p})_{[1,1]}](0)[(E_{p,k_1+q})_{[1,1]}(n_1,n_2)+(E_{q,k_1+p})_{[1,1]}(n_2,n_1)]\\&=&-(1+\dlt_{p,q})[(E_{p,p})_{[1,1]}(n_1,n_2)-(E_{k_1+p,k_1+p})_{[1,1]}(n_2,n_1)]\\& &-(1+\dlt_{p,q})[(E_{q,q})_{[1,1]}(n_2,n_1)-(E_{k_1+q,k_1+q})_{[1,1]}(n_1,n_2)]\in{\cal I}\hspace{3.6cm}(3.92)\end{eqnarray*}
by (2.31) and (2.62). So we have $\bar{R}\bigcap {\cal I}\neq\{0\}$.
Thus (3.91) holds again by Theorem 3.1. We can similarly prove (3.91) if the third element in (3.89) is in ${\cal I}$.

Note that (3.91) implies
$$w=\sum_{j=1}^{k_1}[(E_{j,j})_{[1,1]}(0,1)-(E_{k_1+j,k_1+j})_{[1,1]}(1,0)]\in {\cal I},\eqno(3.93)$$
$$w'=\sum_{j=1}^{k_1}[(E_{j,j})_{[1,1]}(1,0)-(E_{k_1+j,k_1+j})_{[1,1]}(0,1)]\in{\cal I}.\eqno(3.94)$$
Moreover, 
\begin{eqnarray*}& &w(0)[(E_{j_1,k_1+j_2})_{[1,1]}(m,n)+(E_{j_2,k_1+j_1})_{[1,1]}(n,m)]\\&=&-(m+n+2)[(E_{j_1,k_1+j_2})_{[1,1]}(m,n)+(E_{j_2,k_1+j_1})_{[1,1]}(n,m)]\in{\cal I}\hspace{3cm}(3.95)\end{eqnarray*}
\begin{eqnarray*}& &w'(0)[(E_{k_1+j_1,k_1})_{[1,1]}(m,n)+(E_{k_1+j_2,j_1})_{[1,1]}(n,m)]\\&=&(m+n+2)[(E_{k_1+j_1,j_2})_{[1,1]}(m,n)+(E_{k_1+j_2,j_1})_{[1,1]}(n,m)]\in{\cal I}\hspace{3.3cm}(3.96)\end{eqnarray*}
for $j_1,j_2\in\ol{1,k_1}$ and $m,n\in\Bbb{N}$. Thus ${\cal I}=R^{\dg}_{k\times k,1}$. So $R^{\dg}_{k\times k,1}$ is simple.

Assume $k_1>1$. Then $\bar{R}$ is generated by
\begin{eqnarray*}\bar{R}^{(2)}&=&\mbox{span}\{(E_{j_1,j_2})_{[1,1]}(\es,1-\es)-(E_{k_1+j_2,k_1+j_1})_{[1,1]}(1-\es,\es)\\& &\mid j_1,j_2\in\ol{1,k_1},\;\es=0,1\}\subset (R^{\dg}_{k\times k,1})^{(2)}\hspace{6.5cm}(3.97)\end{eqnarray*}
by Theorem 3.1. For any $p,q\in\ol{1,k_1}$ and $m,n\in\Bbb{N}$, we have
\begin{eqnarray*}& &[(E_{k_1+p,q})_{[1,1]}(0,1)+(E_{k_1+q,p})_{[1,1]}(1,0)](0)[(E_{q,q})_{[1,1]}(m,n)-(E_{k_1+q,k_1+q})_{[1,1]}(n,m)]\\&=&(-m-1+m\dlt_{p,q})[(E_{k_1+p,q})_{[1,1]}(m,n)+(E_{k_1+q,p})_{[1,1]}(n,m)],\hspace{3.3cm}(3.98)\end{eqnarray*}
\begin{eqnarray*}& &[(E_{p,k_1+q})_{[1,1]}(0,1)+(E_{q,k_1+p})_{[1,1]}(1,0)](0)[(E_{q,q})_{[1,1]}(m,n)-(E_{k_1+q,k_1+q})_{[1,1]}(n,m)]\\&=&(n+1-n\dlt_{p,q})[(E_{p,k_1+q})_{[1,1]}(m,n)+(E_{k_1+q,p})_{[1,1]}(n,m)]\hspace{4cm}(3.99)\end{eqnarray*}
by (2.31) and (2.62). Thus $R^{\dg}_{k\times k,1}$ is generated by $(R^{\dg}_{k\times k,1})^{(2)}$.

Suppose $k_1=1$. Let $R'$ be the subalgebra of $(R^{\dg}_{k\times k,1})^{(2)}$ generated by (3.80). Then
\begin{eqnarray*}& &[(E_{1,2})_{[1,1]}(2,0)+(E_{1,2})_{[1,1]}(0,2)](0)[(E_{2,1})_{[1,1]}(2,0)+(E_{2,1})_{[1,1]}(0,2)]\\&=&[(E_{1,1})_{[1,1]}(0,2)-(E_{2,2})_{[1,1]}(2,0)]\\& &+7[(E_{1,1})_{[1,1]}(2,0)-(E_{2,2})_{[1,1]}(0,2)]\in R',\hspace{6.3cm}(3.100)\end{eqnarray*}
\begin{eqnarray*}& &[(E_{2,1})_{[1,1]}(2,0)+(E_{2,1})_{[1,1]}(0,2)](0)[(E_{1,2})_{[1,1]}(2,0)+(E_{1,2})_{[1,1]}(0,2)]\\&=&-7[(E_{1,1})_{[1,1]}(0,2)-(E_{2,2})_{[1,1]}(2,0)]\\& &-[(E_{1,1})_{[1,1]}(2,0)-(E_{2,2})_{[1,1]}(0,2)]\in R'\hspace{6.6cm}(3.101)\end{eqnarray*}
by (2.31) and (2.62). Solving (3.100) and (3.101), we obtain
$$(E_{1,1})_{[1,1]}(2,0)-(E_{2,2})_{[1,1]}(0,2)\in R'.\eqno(3.102)$$
By Theorem 3.1, $\bar{R}\subset R'$. Moreover,
\begin{eqnarray*}& &[(E_{1,1})_{[1,1]}(0,1)-(E_{2,2})_{[1,1]}(1,0)](1)[(E_{2,1})_{[1,1]}(2,0)+(E_{2,1})_{[1,1]}(0,2)]\\&=&-[(E_{2,1})_{[1,1]}(1,0)+(E_{2,1})_{[1,1]}(0,1)]\in R',\hspace{6.5cm}(3.103)\end{eqnarray*}
\begin{eqnarray*}& &[(E_{1,1})_{[1,1]}(0,1)-(E_{2,2})_{[1,1]}(1,0)](1)[(E_{1,2})_{[1,1]}(2,0)+(E_{1,2})_{[1,1]}(0,2)]\\&=&-3[(E_{1,2})_{[1,1]}(1,0)+(E_{1,2})_{[1,1]}(0,1)]\in R'\hspace{6.5cm}(3.104)\end{eqnarray*}
by (2.31) and (2.62). Thus $R'=R^{\dg}_{k\times k,1}$ (3.98), (3.99), (3.103) and (3.104).

Next we consider $R^{\dg}_{k\times k,2}$.
Let ${\cal I}$ be a nonzero ideal of $R_{k\times k,2}^{\dg}$. For $p,q,r\in\ol{1,k_1}$ and $m,n_1,n_2\in\Bbb{N}$, we have
\begin{eqnarray*}& &[E_{r,r}(0,m)+E_{k_1+r,k_1+r}(m,0)](0)[E_{p,q}(n_1,n_2)+E_{k_1+q,k_1+p}(n_2,n_1)]\\&=&[\dlt_{r,p}(n_1+1)(^{-n_1-2}_{\;\;\;\;m})+\dlt_{r,q}(n_2+1)(^{n_2}_m)][E_{p,q}(n_1,n_2)+
E_{k_1+q,k_1+p}(n_2,n_1)],\hspace{1cm}(3.105)\end{eqnarray*}
\begin{eqnarray*}& &[E_{r,r}(0,m)+E_{k_1+r,k_1+r}(m,0)](0)[E_{p,k_1+q}(n_1,n_2)-E_{q,k_1+p}(n_2,n_1)]\\&=&[\dlt_{r,p}(n_1+1)(^{-n_1-2}_{\;\;\;\;m})+\dlt_{r,q}(n_2+1)(^{-n_2-2}_{\;\;\;\;m})]\\& &[E_{p,k_1+q}(n_1,n_2)-E_{q,k_1+p}(n_2,n_1)],\hspace{7.6cm}(3.106)\end{eqnarray*}
\begin{eqnarray*}& &[E_{r,r}(0,m)+E_{k_1+r,k_1+r}(m,0)](0)[E_{k_1+p,q}(n_1,n_2)-E_{k_1+q,p}(n_2,n_1)]\\&=&[\dlt_{r,p}(n_1+1)(^{n_1}_m)+\dlt_{r,q}(n_2+1)(^{n_2}_m)][E_{k_1+p,q}(n_1,n_2)-E_{k_1+q,p}(n_2,n_1)]\hspace{1.7cm}(3.107)\end{eqnarray*}
by (2.31) and (2.62). Moreover,
$$\dlt_{r,p}(n_1+1)(^{-n_1-2}_{\;\;\;\;m})+\dlt_{r,q}(n_2+1)(^{n_2}_m)={1\over m!}{d^{m+1}\over dx^{m+1}}(-\dlt_{r,p}x^{-n_1-1}+\dlt_{r,q}x^{n_2+1})|_{x=1},\eqno(3.108)$$
\begin{eqnarray*}& &\dlt_{r,p}(n_1+1)(^{-n_1-2}_{\;\;\;\;m})+\dlt_{r,q}(n_2+1)(^{-n_2-2}_{\;\;\;\;m})\\&=&-{1\over m!}{d^{m+1}\over dx^{m+1}}(\dlt_{r,p}x^{-n_1-1}+\dlt_{r,q}x^{-n_2-1})|_{x=1},\hspace{6.4cm}(3.109)\end{eqnarray*}
$$\dlt_{r,p}(n_1+1)(^{n_1}_m)+\dlt_{r,q}(n_2+1)(^{n_2}_m)={1\over m!}{d^{m+1}\over dx^{m+1}}(\dlt_{r,p}x^{n_1+1}+\dlt_{r,q}x^{n_2+1})|_{x=1}.\eqno(3.110)$$

Let ${\cal I}$ be a nonzero ideal of $R^{\dg}_{k\times k,2}$. By Lemma 2.1, (3.105)-(3.110) and the Taylor's theorem at $x=1$ in calculus (cf. (3.9)-(3.12)),  ${\cal I}$ contains at least one of the following elements:
$$E_{p,q}(n_1,n_2)+E_{k_1+q,k_1+p}(n_2,n_1),\;\;E_{p,k_1+q}(n_1,n_2)-E_{q,k_1+p}(n_2,n_1),\eqno(3.111)$$
$$E_{k_1+p,q}(n_1,n_2)-E_{k_1+q,p}(n_2,n_1)\eqno(3.112)$$
for some $p,q\in\ol{1,k_1}$ and $n_1,n_2\in\Bbb{N}$, where $p\neq q$ or $n_1\neq n_2$ in the second and third elements.
 Note that the subspace
$$\hat{R}=\mbox{span}\{E_{j_1,j_2}(m,n)+E_{k_1+j_2,k_1+j_1}(n,m)\mid j_1,j_2\in\ol{1,k_1},\;m,n\in\Bbb{N}\}\eqno(3.113)$$
forms a subalgebra of $R^{\dg}_{k\times k,2}$ that is isomorphic to $R_{k_1\times k_1,2}$. By Theorem 3.1, $R_{k_1\times k_1,2}$ is simple. Hence
$$\hat{R}\in{\cal I}\eqno(3.114)$$
if the first element in (3.111) is in $I$. Assume that the second element in (3.111) is in ${\cal I}$. Then we have
\begin{eqnarray*}& &[E_{k_1+p,q}(0,1)-E_{k_1+q,p}(1,0)](0)[E_{p,k_1+q}(n_1,n_2)-E_{q,k_1+p}(n_2,n_1)]\\&=&(n_2+1)(n_2+2+n_2\dlt_{p,q})[E_{p,p}(n_1,n_2)+E_{k_1+p,k_1+p}(n_2,n_1)]\\& &-(n_1+1)(n_1+(n_1+2)\dlt_{p,q})[E_{q,q}(n_2,n_1)+E_{k_1+q,k_1+q}(n_1,n_2)]\in{\cal I}\hspace{2cm}(3.115)\end{eqnarray*}
by (2.31) and (2.62). When $p=q$, the coefficient of the lower term is
$$2(n_2+1)^2-2(n_1+1)^2=2(n_2-n_1)(n_1+n_2+2),\eqno(3.116)$$
which is zero only if $n_1=n_2$.
So we have $\hat{R}\bigcap {\cal I}\neq\{0\}$.
Thus (3.113) holds again by Theorem 3.1. We can similarly prove (3.113) if the element in (3.112) is in ${\cal I}$.
Note that $I_k\in \hat{R}$ (cf. (3.5)). Hence ${\cal I}=R^{\dg}_{k\times k,2}$ by (3.19) with $j=\ell=0$. Therefore $R^{\dg}_{k\times k,2}$ is simple.

Assume $k_1>1$. Then $\hat{R}$ is generated by
$$\hat{R}^{(2)}=\mbox{span}\{E_{j_1,j_2}+E_{k_1+j_2,k_1+j_1}\mid j_1,j_2\in\ol{1,k_1}\}\subset (R^{\dg}_{k\times k,2})^{(2)}\eqno(3.117)$$
by Theorem 3.1 (cf. (3.5)). For $m,n\in\Bbb{N}$ and $p,q\in\ol{1,k_1}$ such that $p\neq q$, we have
\begin{eqnarray*}& &[(E_{p,p}+E_{k_1+p,k_1+p})(-1)]^m[(E_{q,q}+E_{k_1+q,k_1+q})(-1)]^n(E_{p,k_1+q}-E_{q,k_1+p})\\&=&m!n!(E_{p,k_1+q}(m,n)-E_{q,k_1+p}(n,m)),\hspace{7.2cm}(3.118)\end{eqnarray*}
\begin{eqnarray*}& &[(E_{p,p}+E_{k_1+p,k_1+p})(-1)]^m[(E_{q,q}+E_{k_1+q,k_1+q})(-1)]^n(E_{k_1+p,q}-E_{k_1+q,p})\\&=&m!n!(E_{k_1+p,q}(m,n)-E_{k_1+q,p}(n,m)),\hspace{7.2cm}(3.119)\end{eqnarray*}
\begin{eqnarray*}& &(E_{p,q}+E_{k_1+q,k_1+p})(0)(E_{p,k_1+q}(m,n)-E_{q,k_1+p}(n,m))\\&=&(n+1)(E_{p,k_1+p}(m,n)-E_{p,k_1+p}(n,m)).\hspace{6.8cm}(3.120)\end{eqnarray*}
\begin{eqnarray*}& &(E_{p,q}+E_{k_1+q,k_1+p})(0)(E_{k_1+q,p}(m,n)-E_{k_1+p,q}(n,m))\\&=&(n+1)(E_{k_1+q,q}(m,n)-E_{k_1+q,q}(n,m))\hspace{6.8cm}(3.121)\end{eqnarray*}
by (2.31) and (2.62). Thus $R^{\dg}_{k\times k,2}$ is generated by $(R^{\dg}_{k\times k,2})^{(2)}$.

Next we consider the case when $k_1=1$. Let $R'$ be the subalgebra of $R^{\dg}_{k\times k,2}$ generated by (3.82).
 Note 
\begin{eqnarray*}& &(E_{1,2}(1,0)-E_{1,2}(0,1))(0)(E_{2,1}(0,1)-E_{2,1}(1,0))\\&=&2(E_{1,1}(0,1)+E_{2,2}(1,0))-8(E_{1,1}(1,0)+E_{2,2}(0,1))\in R',\hspace{3.5cm}(3.122)\end{eqnarray*}
$$(I_2)(-1)(I_2)=E_{1,1}(1,0)+E_{1,1}(0,1)+E_{2,2}(0,1)+E_{2,2}(1,0)\in R'.\eqno(3.123)$$
Thus
$$E_{1,1}(1,0)+E_{2,2}(0,1),\; E_{1,1}(0,1)+E_{2,2}(1,0)\in R'.\eqno(3.124)$$
By Theorem 3.1,
$$\hat{R}\subset R'.\eqno(3.125)$$
For $m,n\in\Bbb{N}$, we have
\begin{eqnarray*}& &(E_{1,2}(1,0)-E_{1,2}(0,1))(0)(E_{1,1}(m,n)+E_{2,2}(n,m))\\&=&2(n+1)^2[E_{1,2}(n,m)-E_{1,2}(m,n)]\hspace{7.7cm}(3.126)\end{eqnarray*}
\begin{eqnarray*}& &(E_{2,1}(1,0)-E_{2,1}(0,1))(0)(E_{1,1}(m,n)+E_{2,2}(n,m))\\&=&2(m+1)^2[E_{2,1}(m,n)-E_{2,1}(n,m)]\hspace{7.6cm}(3.127)\end{eqnarray*}
by (2.31) and (2.62). Thus $R'=R^{\dg}_{k\times k,2}.\qquad\Box$
\psp

{\bf Remark 3.4}. (a) The algebra $R_{1\times 1,2}$ is the well-known $W_{\infty}$ algebra without center (cf. [Ba]) and  the algebra $R_{1\times 1,1}$ is the well-known $W_{1+\infty}$ algebra without center (cf. [PRS]) in mathematical physics. The more general algebra $R_{k\times k,1}$ is the $W_{1+\infty}(gl_k)$ studied by van de Leur [V] without center related to $k$ component KP hierarchy. 

(b)  Let $(R,\ptl,Y^+(\cdot,z))$ be a $\G$-weighted conformal algebra (cf. (1.7), (1.8)). For each $\al\in\G$, we define the homogeneous algebraic operation $\odot$ on $R^{(\al)}$ by
$$u\odot v= u(0)v \qquad\for\;\;u,v\in R^{(\al)}\eqno(3.128)$$
(cf. (1.13)). 

The homogeneous subalgebraic structure $((R_{k\times k,2\ell+2})^{(2\ell+2)},\odot)$ of the algebra $R_{k\times k,2\ell+2}$ with $\ell\in\Bbb{N}$ has the property:
$$u(0,2\ell)\odot v(0,2\ell)=(2\ell+1)(uv+vu)(0,2\ell)\qquad\for\;\;u,v\in M_{k\times k}(\Bbb{F}).\eqno(3.129)$$
So $((R_{k\times k,2\ell+2})^{(2\ell+2)},\odot)$ is isomorphic to the simple Jordan algebra of type $A_k$. Thus by the generator property in Theorem 3.1, the simple conformal algebra $R_{k\times k,2\ell+2}$ with $\ell\in\Bbb{N}$ and $k>1$ can be viewed to be generated by the simple Jordan algebra $((R_{k\times k,2\ell+2})^{(2\ell+2)},\odot)$ of type $A_k$.

The homogeneous subalgebraic structure $((R_{k\times k,2\ell+3})^{(2\ell+3)},\odot)$ of the algebra $R_{k\times k,2\ell+3}$ with $\ell\in\Bbb{N}$ has the property:
$$u(0,2\ell+1)\odot v(0,2\ell+1)=(2\ell+2)(-uv+vu)(0,2\ell)\qquad\for\;\;u,v\in M_{k\times k}(\Bbb{F}).\eqno(3.130)$$
So $((R_{k\times k,2\ell+3})^{(2\ell+3)},\odot)$ is isomorphic to the Lie algebra $gl_k(\Bbb{F})$. Thus by the generator property in Theorem 3.1, the simple conformal algebra $R_{k\times k,2\ell+3}$ with $\ell\in\Bbb{N}$ and $k>1$ are generated by the Lie algebra $((R_{k\times k,2\ell+3})^{(2\ell+3)},\odot)$ of type $gl_k(\Bbb{F})$.

Similarly, we can view that the simple conformal algebras $R^{\ast}_{k\times k,2}$ with $k>1$ are generated by the simple Jordan algebra $((R^{\ast}_{k\times k,2})^{(2)},\odot)$ of type $B_k$ and $R^{\dg}_{k\times k,2}$  with even $k>1$ are generated by the simple Jordan algebra $((R^{\dg}_{k\times k,2})^{(2)},\odot)$ of type $C_{k/2}$.

\section{Simple Conformal Superalgebras}

In this section, we shall construct three families of conformal superalgebras of finite growth with nonzero odd part from  matrix algebras and prove their simplicity.

Let us go back to the general construction of $(R({\cal A}),Y^+(\cdot,z))$ in (2.23)-(2.33). We let ${\cal A}=M_{k\times k}(\Bbb{F})$ the $k\times k$ matrix algebra. Assume that 
$$k=k_1+k_2\eqno(4.1)$$
for some fixed positive integers $k_1$ and $k_2$. Then ${\cal A}$ has the following $\Bbb{Z}_2$-grading:
$${\cal A}_0=\left\{\left(\begin{array}{cc}A_{1,1}&\\ &A_{2,2}\end{array}\right)\mid A_{1,1}\in M_{k_1\times k_1}(\Bbb{F}),\;A_{2,2}\in M_{k_2\times k_2}(\Bbb{F})\right\},\eqno(4.2)$$
$${\cal A}_1=\left\{\left(\begin{array}{cc}&A_{1,2}\\ A_{2,1}&\end{array}
\right)\mid A_{1,2}\in M_{k_1\times k_2}(\Bbb{F}),\;A_{2,1}\in M_{k_2\times k_1}(\Bbb{F})\right\}.\eqno(4.3)$$
For $\ell\in\Bbb{N}$, we define 
$$R_{[k_1,k_2],\ell+1}=(R_{[k_1,k_2],\ell+1})_0+(R_{[k_1,k_2],\ell+1})_1\subset R({\cal A})\eqno(4.4)$$
by
\begin{eqnarray*}(R_{[k_1,k_2],\ell+1})_0&=&\mbox{span}\:\{E_{p_1,q_1}(m,n+\ell),(E_{k_1+p_2,k_1+q_2})_{[1,1]}(m,\ell+n)\\& &\mid p_1,q_1\in\ol{1,k_1},\;p_2,q_2\in\ol{1,k_2},\;m,n\in\Bbb{N}\},\hspace{4.6cm}(4.5)\end{eqnarray*}
\begin{eqnarray*}(R_{[k_1,k_2],\ell+1})_1&=&\mbox{span}\:\{(E_{p_1,k_1+p_2})_{[0,1]}(m,\ell+n),\;(E_{k_1+p_2,p_1})_{[1,0]}(m,\ell+n)\\& &\mid p_1\in\ol{1,k_1},\;p_2\in\ol{1,k_2},\;m,n\in\Bbb{N}\}.\hspace{5.8cm}(4.6)\end{eqnarray*}
Then each $R_{[k_1,k_2],\ell+1}$ forms a conformal sub-superalgebra of $R({\cal A})$ for $\ell\in\Bbb{N}$.
\psp

{\bf Theorem 4.1}. {\it Let} $\ell\in\Bbb{N}$. {\it The algebra} $(R_{[k_1,k_2],\ell+1},\ptl,Y^+(\cdot,z))$ {\it is simple. When} $k>2$, {\it  it is generated by} 
\begin{eqnarray*}R_{[k_1,k_2],\ell+1}^{(\ell+3/2)}&=&\mbox{span}\;\{(E_{p_1,k_1+p_2})_{[0,1]}(0,\ell),(E_{k_1+p_2,p_1})_{[1,0]}(0,\ell)\\& &\mid p_1\in\ol{1,k_1},\;p_2\in\ol{1,k_2}\}.\hspace{8.1cm}(4.7)\end{eqnarray*}
{\it If} $k=2$, {\it then it is is generated by} 
$$\{(E_{1,2})_{[0,1]}(0,\ell),(E_{2,1})_{[1,0]}(0,\ell),(E_{2,2})_{[1,1]}(0,\ell+1)\}.\eqno(4.8)$$

{\it Proof}. We fix $\ell\in\Bbb{N}$. We let 
$$R_{0,0}=\mbox{span}\:\{E_{p_1,q_1}(m,n+\ell)\mid p_1,q_1\in\ol{1,k_1},\;m,n\in\Bbb{N}\},\eqno(4.9)$$
$$R_{0,1}=\mbox{span}\:\{(E_{p_2,q_2})_{[1,1]}(m,n+\ell)\mid p_2,q_2\in\ol{1,k_1},\;m,n\in\Bbb{N}\},\eqno(4.10)$$
Then $R_{0,0}$ and $R_{0,1}$ form conformal subalgebras of $R_{[k_1,k_2],\ell+1}$. In fact,
$$R_{0,0}\cong R_{k_1\times k_1,\ell+2},\;\;R_{0,1}\cong R_{k_2\times k_2,\ell+1}\eqno(4.11)$$
by (3.32). Hence $R_{0,0}$ and $R_{0,1}$ both are simple subalgebras. Moreover,  we have
$$Y^+(u_{[i,i]}(m_1,m_2),z)v_{[j,j]}(n_1,n_2)=0\eqno(4.12)$$
for $u,v\in {\cal A}_0,\;m_1,m_2,n_1,n_2\in\Bbb{N},\;i,j\in\Bbb{Z}_2,i\neq j$ by (2.31).  

For $p_1,r_1\in\ol{1,k_1},\;p_2,r_2\in\ol{1,k_2}$ and $m,n_1,n_2\in\Bbb{N}$ with $m,n_2\geq \ell$, we obtain
\begin{eqnarray*}& &(E_{r_1,r_1}(0,m))(0)(E_{p,q})_{[i,j]}(n_1,n_2)\\&=&
[\dlt_{r_1,p}\dlt_{i,0}(n_1+1)(^{-n_1-2}_{\;\;\;\;m})+\dlt_{r_1,q}\dlt_{j,0}(n_2+1)(^{n_2}_m)](E_{p,q})_{[i,j]}(n_1,n_2),\hspace{2.8cm}(4.13)\end{eqnarray*}
\begin{eqnarray*}& &[(E_{k_1+r_2,k_1+r_2})_{[1,1]}(0,m)](0)(E_{p,q})_{[i,j]}(n_1,n_2)\\&=&[\dlt_{k_1+r_2,p}\dlt_{i,1}(^{-n_1-1}_{\;\;\;\;m})-\dlt_{k_1+r_2,q}\dlt_{j,1}(^{n_2}_m)](E_{p,q})_{[i,j]}(n_1,n_2)\hspace{4.8cm}(4.14)\end{eqnarray*}
by (2.31), where $(p,q)=(p_1,k_1+p_2), [i,j]=[0,1]$ or $(p,q)=(k_1+p_2,p_1),[i,j]=[1,0]$. Note that
\begin{eqnarray*}& &
\dlt_{r_1,p}\dlt_{i,0}(n_1+1)(^{-n_1-2}_{\;\;\;\;m})+\dlt_{r_1,q}\dlt_{j,0}(n_2+1)(^{n_2}_m)\\&=&{1\over m!}{d^{m+1}\over dx^{m+1}}(\dlt_{r_1,q}\dlt_{j,0}x^{n_2+1}-\dlt_{r_1,p}\dlt_{i,0}x^{-n_1-1})|_{x=1},\hspace{5.8cm}(4.15)\end{eqnarray*}
\begin{eqnarray*}& &\dlt_{k_1+r_2,p}\dlt_{i,1}(^{-n_1-1}_{\;\;\;\;m})-\dlt_{k_1+r_2,q}\dlt_{j,1}(^{n_2}_m)\\&=&{1\over m!}{d^m\over dx^m}(\dlt_{k_1+r_2,p}\dlt_{i,1}x^{-n_1-1}-\dlt_{k_1+r_2,q}\dlt_{j,1}x^{n_2})|_{x=1}.\hspace{5.6cm}(4.16)\end{eqnarray*}

Let ${\cal I}$ be a nonzero ideal of $R_{[k_1,k_2],\ell+1}$.
By Lemma 2.1, (4.12)-(4.16) and Taylor's Theorem at $x=1$ in calculus, we have $R_{0,0}\subset {\cal I}$ or $R_{0,1}\subset {\cal I}$ or 
$$(E_{p,k_1+q})_{[0,1]}(n_1,n_2)\in{\cal I}\;\mbox{or}\;(E_{k_1+q,p})_{[1,0]}(n_1,n_2)\in{\cal I}\eqno(4.17)$$
for some $p\in\ol{1,k_1},\;q\in\ol{1,k_2}$ and $n_1,n_2\in\Bbb{N}$ with $n_2\geq\ell$. Note that either of the first two cases implies (4.17) by (4.13) and (4.14). Without loss of generality, we can assume the first case in (4.17). Furthermore, 
\begin{eqnarray*}& &[(E_{k_1+q,p})_{[1,0]}(0,\ell)](-1/2)(E_{p,k_1+q})_{[0,1]}(n_1,n_2)\\&=&(^{n_2}_{\:\ell})E_{p,p}(n_1,n_2)+(n_1+1)(^{-n_1-2}_{\:\;\;\;\:\ell})(E_{k_1+q,k_1+q})_{[1,1]}(n_1+1,n_2)\in{\cal I}\hspace{2.6cm}(4.18)\end{eqnarray*}
by (2.31) and (2.62). Hence ${\cal I}\bigcap R_{0,0}\neq \{0\}$ and
${\cal I}\bigcap R_{0,1}\neq \{0\}$ by (4.12). Since $R_{0,0}$ and $R_{0,1}$ are simple, we have $R_{0,0},R_{0,1}\subset{\cal I}$. Hence ${\cal I}=R_{[k_1,k_2],\ell}$ by (4.13) and (4.14). Thus $R_{[k_1,k_2],\ell+1}$ is simple.

Let $R'$ be the subalgebra of $(R_{[k_1,k_2],\ell+1})^{(\ell+3/2)}$ when $k>2$ and by (4.8) when $k=2$. Note that for $p_1,p_2\in\ol{1,k_1}$ and $q_1,q_2\in\ol{1,k_2}$, we have
\begin{eqnarray*}& &[(E_{p_1,k_1+1})_{[0,1]}(0,\ell)](-1/2)[(E_{k_1+1,p_2})_{[1,0]}(0,\ell)]\\&=&(-1)^{\ell}E_{p_1,p_2}(0,\ell)-(\ell+1)^2\dlt_{p_1,p_2}(E_{k_1+1,k_1+1})_{[1,1]}(0,\ell+1),\hspace{3.9cm}(4.19)\end{eqnarray*}
\begin{eqnarray*}& &[(E_{1,k_1+q_2})_{[0,1]}(0,\ell)](1/2)[(E_{k_1+q_1,1})_{[1,0]}(0,\ell)]\\&=&-(\ell+1)(E_{k_1+q_1,k_1+q_2})_{[1,1]}(0,\ell),\hspace{8.4cm}(4.20)\end{eqnarray*}
\begin{eqnarray*}& &[(E_{1,k_1+p_2})_{[0,1]}(0,\ell)](-1/2)[(E_{k_1+q_1,1})_{[1,0]}(0,\ell)]\\&=&(-1)^{\ell}\dlt_{q_1,q_2}E_{1,1}(0,\ell)-(\ell+1)^2(E_{k_1+q_1,k_1+q_2})_{[1,1]}(0,\ell+1)\hspace{4.1cm}(4.21)\end{eqnarray*}  
by (2.31) and (2.62).

If $k_1>1$, then
$$E_{1,1}(0,\ell)-E_{2,2}(0,\ell)\in R'\eqno(4.22)$$
by (4.19). Moreover,
$$[(E_{1,1})(0,\ell)-E_{2,2}(0,\ell)](-1)(E_{k_1+1,1})_{[1,0]}(0,\ell)=(\ell+1)^2(E_{k_1+1,1})_{[1,0]}(0,\ell+1)\eqno(4.23)$$
by (2.31) and (2.62). So $(E_{k_1+1,1})_{[1,0]}(0,\ell+1)\in R'$. Moreover,
\begin{eqnarray*}& &[(E_{1,k_1+1})_{[0,1]}(0,\ell)](1/2)(E_{k_1+1,1})_{[1,0]}(0,\ell+1)\\&=&-(\ell+1)(\ell+2)(E_{k_1+1,k_1+1})_{[1,1]}(0,\ell+1)\hspace{6.9cm}(4.24)\end{eqnarray*}
by (2.31) and (2.62). Hence
$$(E_{k_1+1,k_1+1})_{[1,1]}(0,\ell+1)\in R',\eqno(4.25)$$
which implies
$$E_{p_1,p_2}(0,\ell)\in R'\qquad\for\;\;p_1,p_2\in \ol{1,k_1}\eqno(4.26)$$
by (4.19). Furthermore,  Theorem 3.1 and (4.26) implies $R_{0,0}\in R'$. For $p\in\ol{1,k_1},\;q\in\ol{1,k_2}$ and $m,n\in\Bbb{N}$ with $n\geq \ell$, we have
$$[(E_{k_1+q,p})_{[1,0]}(0,\ell)](1/2)E_{p,p}(m,n)=(m+1)(^{-m-1}_{\;\;\;\;\ell})(E_{k_1+q,p})_{[1,0]}(m,n),\eqno(4.27)$$
$$[(E_{p,k_1+q})_{[0,1]}(0,\ell)](1/2)E_{p,p}(m,n)=(n+1)(^n_{\ell})(E_{p,k_1+q})_{[0,1]}(m,n)\eqno(4.28)$$
by (2.31) and (2.62). Thus
$$(R_{[k_1,k_2],\ell+1})_1\subset R'\eqno(4.29)$$
(cf. (4.6)). For $q_1,q_2\in\ol{1,k_2}$ and $m,n\in\Bbb{N}$ with $n\geq \ell$, we have
\begin{eqnarray*}& &[(E_{k_1+q_1,1})_{[1,0]}(0,\ell)](1/2)(E_{1,k_1+q_2})_{[0,1]}(m,n)\\&=&(m+1)(^{-m-2}_{\;\;\;\;\:\ell})(E_{k_1+q_1,k_1+q_2})_{[1,1]}(m,n)+(^{n-1}_{\;\;\;\ell})E_{1,1}(m,n-1)\in R'\hspace{2.4cm}(4.30)\end{eqnarray*}
by (2.31) and (2.62). Since $E_{1,1}(m,n-1)\in R'$, we have
$$(E_{k_1+q_1,k_1+q_2})_{[1,1]}(m,n)\in R'\qquad\for\;\;q_1,q_2\in\ol{1,k_2},\;m,n\in\Bbb{N},\;n\geq \ell.\eqno(4.31)$$
Therefore $R'=R_{[k_1,k_2],\ell+1}$. 

Assume $k_2>1$. If $\ell>0$, then we get (4.31) by Theorem 3.1,  (3.32) and (4.20). If $\ell=0$, we have
$$(E_{k_1+1,k_1+2})_{[1,1]}(0,1),(E_{k_1+2,k_1+1})_{[1,1]}(0,1)\in R',\eqno(4.32)$$
$$(E_{k_1+1,k_1+1})_{[1,1]}(0,1)+(E_{k_1+2,k_1+2})_{[1,1]}(0,1)-2E_{1,1}\in R'\eqno(4.33)$$
by (3.5) and (4.21). Moreover,
\begin{eqnarray*}& &[(E_{k_1+1,k_1+2})_{[1,1]}(0,1)](0)(E_{k_1+2,k_1+1})_{[1,1]}(0,1)\\&=&-[(E_{k_1+1,k_1+1})_{[1,1]}(0,1)+(E_{k_1+2,k_1+2})_{[1,1]}(0,1)]\in R'\hspace{4.7cm}(4.34)\end{eqnarray*}
by (2.31) and (2.62). Thus $E_{1,1}\in R'$ by (4.33) and (4.34). Hence
$$(E_{k_1+q_1,k_1+q_2})_{[1,1]}(0,1)\in R'\qquad\for\;\;q_1,q_2\in\ol{1,k_2}\eqno(4.35)$$
by (4.21). Symmetrically, we can prove $(E_{k_1+q_1,k_1+q_2})_{[1,1]}(1,0)\in R'$ for $q_1,q_2\in\ol{1,k_2}$. Expression (4.31) holds again by Theorem 3.1.
 For $p\in\ol{1,k_1},\;q\in\ol{1,k_2}$ and $m,n\in\Bbb{N}$ with $n\geq \ell$, we have
$$[(E_{k_1+q,p})_{[1,0]}(0,\ell)](-1/2)(E_{k_1+q,k_1+q})_{[1,1]}(m,n)=-(^n_{\ell})(E_{k_1+q,p})_{[1,0]}(m,n),\eqno(4.36)$$
$$[(E_{p,k_1+q})_{[0,1]}(0,\ell)](-1/2)(E_{k_1+q,k_1+q})_{[1,1]}(m,n)=(^{-m-1}_{\;\;\;\;\ell})(E_{p,k_1+q})_{[1,0]}(m,n),\eqno(4.37)$$
by (2.31) and (2.62). Thus (4.29) holds. For $p_1,p_2\in\ol{1,k_1}$ and $m,n\in\Bbb{N}$ with $n\geq \ell$, we have
\begin{eqnarray*}& &[(E_{k_1+1,p_2})_{[1,0]}(0,\ell)](-1/2)(E_{p_1,k_1+1})_{[0,1]}(m,n)\\&=&(^n_{\ell})E_{p_1,p_2}(m,n)+(m+1)(^{-m-2}_{\;\;\;\;\ell})\dlt_{p_1,p_2}(E_{k_1+1,k_1+1})_{[1,1]}(m+1,n)\in R'\hspace{1.9cm}(4.38)\end{eqnarray*}
by (2.31) and (2.62). Since $(E_{k_1+1,k_1+1})_{[1,1]}(m+1,n)\in R'$, we have
$$E_{p_1,p_2}(m,n)\in R'\qquad\for\;\;p_1,p_2\in\ol{1,k_1},\;m,n\in\Bbb{N},\;n\geq \ell.\eqno(4.39)$$
Therefore $R'=R_{[k_1,k_2],\ell+1}$.

Let $k_1=k_2=1$. By (4.8) and  (4.19),
$$E_{1,1}(0,\ell)\in R'.\eqno(4.40)$$
Moreover,
$$(E_{1,1}(0,\ell))(-1)(E_{2,1})_{[1,0]}(0,\ell)=(\ell+1)^2(E_{2,1})_{[1,0]}(0,\ell+1)\in R'\eqno(4.41)$$
by (2.31) and (2.62). Furthermore,
\begin{eqnarray*}& &(E_{1,1}(0,\ell))(0)[(E_{1,2})_{[0,1]}(0,\ell)](-1/2)(E_{2,1})_{[1,0]}(0,\ell+1)\\&=&(E_{1,1}(0,\ell))(0)[(-1)^{\ell}E_{1,1}(0,\ell+1)-{(\ell+2)^2(\ell+1)\over 2}(E_{2,2})_{[1,1]}(0,\ell+2)]\\&=& (-1)^{\ell}(\ell+1)((-1)^{\ell}+\ell+2)E_{1,1}(0,\ell+1)\in R'\hspace{5.8cm}(4.42)\end{eqnarray*}
by (4.41). So $E_{1,1}(0,\ell+1)\in R'$. Hence $R_{0,0}\in R'$ Theorem 3.1. Hence $R'=R_{[k_1,k_2],\ell+1}$ (4.27)-(4.31).$\qquad\Box$
\psp

Let $\sgm_1:A\mapsto A^t$ be the transpose map of matrices. Then $\sgm_1$ is an involutive anti-isomorphism of $M_{k\times k}(\Bbb{F})$ preserving the $\Bbb{Z}_2$-grading in (4.2) and (4.3). Thus we have the following subalgebra of $R_{[k_1, k_2],1}$: $R_{[k_1,k_2]}^{\ast}=(R_{[k_1,k_2]}^{\ast})_0+(R_{[k_1,k_2]}^{\ast})_1$ with
\begin{eqnarray*}& &(R_{[k_1,k_2]}^{\ast})_0=\mbox{span}\:\{E_{p_1,q_1}(m,n)+E_{q_1,p_1}(n,m),(E_{k_1+p_2,k_1+q_2})_{[1,1]}(m,n)\\& &-(E_{k_1+q_2,k_1+p_2})_{[1,1]}(n,m)\mid p_1,q_1\in\ol{1,k_1},\;p_2,q_2\in\ol{1,k_2},\;m,n\in\Bbb{N}\},\hspace{2.4cm}(4.43)\end{eqnarray*}
\begin{eqnarray*}& &(R_{[k_1,k_2]}^{\ast})_1=\mbox{span}\:\{(E_{p_1,k_1+p_2})_{[0,1]}(m,n)\\& &+(E_{k_1+p_2,p_1})_{[1,0]}(n,m)\mid p_1\in\ol{1,k_1},\;p_2\in\ol{1,k_2},\;m,n\in\Bbb{N}\}.\hspace{4.1cm}(4.44)\end{eqnarray*}
\psp

{\bf Theorem 4.2}.  {\it The algebra} $(R^{\ast}_{[k_1,k_2]},\ptl,Y^+(\cdot,z))$ {\it is simple. Moreover, it is generated by} 
$$(R_{[k_1, k_2]}^{\ast})^{(3/2)}=\mbox{span}\:\{(E_{p_1,k_1+p_2})_{[0,1]}+(E_{k_1+p_2,p_1})_{[1,0]}\mid p_1\in\ol{1,k_1},\;p_2\in\ol{1,k_2}\}\eqno(4.45)$$
 {\it (cf. (2.61)) if} $k>2$ {\it and by} 
$$\{(E_{1,2})_{[0,1]}+(E_{2,1})_{[1,0]},(E_{2,2})_{[1,1]}(0,1)-(E_{2,2})_{[1,1]}(1,0)\}\eqno(4.46)$$
{\it if} $k=2$.
\psp

{\it Proof}. For $r_1,p\in\ol{1,k_1},\;q\in\ol{1,k_2}$ and $m,n_1,n_2\in\Bbb{N}$, we have
\begin{eqnarray*}& &(E_{r_1,r_1}(0,m)+E_{r_1,r_1}(m,0))(0)[(E_{p,k_1+q})_{[0,1]}(n_1,n_2)+(E_{k_1+q,p})_{[1,0]}(n_2,n_1)]\\&=&\dlt_{r_1,p}(n_1+1)((^{-n_1-2}_{\;\;\;\;m})+(^{n_1}_m))[(E_{p,k_1+q})_{[0,1]}(n_1,n_2)+(E_{k_1+q,p})_{[1,0]}(n_2,n_1)],\hspace{1.2cm}(4.47)\end{eqnarray*}
\begin{eqnarray*}& &(E_{k_1+r_2,k_1+r_2})_{[1,1]}(0,m)-(E_{k_1+r_2,k_1+r_2})_{[1,1]}(m,0))(0)\\& &[(E_{p,k_1+q})_{[0,1]}(n_1,n_2)+(E_{k_1+q,p})_{[1,0]}(n_2,n_1)]\\&=&\dlt_{r_2,q}((^{-n_2-1}_{\;\;\;m})-(^{n_2}_m))[(E_{p,k_1+q})_{[0,1]}(n_1,n_2)+(E_{k_1+q,p})_{[1,0]}(n_2,n_1)],\hspace{2.7cm}(4.48)\end{eqnarray*}
\begin{eqnarray*}& &[(E_{k_1+q,p})_{[1,0]}+(E_{p,k_1+q})_{[0,1]}](-1/2)[(E_{p,k_1+q})_{[0,1]}(n_1,n_2)+(E_{k_1+q,p})_{[1,0]}(n_2,n_1)]\\&=&(n_1+1)[(E_{k_1+q,k_1+q})_{[1,1]}(n_1+1,n_2)-(E_{k_1+q,k_1+q})_{[1,1]}(n_2,n_1+1)]\\& &+E_{p,p}(n_1,n_2)+E_{p,p}(n_2,n_1)\hspace{9.1cm}(4.49)\end{eqnarray*}
(cf. (2.61)) by (2.31) and (2.62). Note that
$$\dlt_{r_1,p}(n_1+1)((^{-n_1-2}_{\;\;\;\;\:m})+(^{n_1}_m))={\dlt_{r_1,p}\over m!}{d^{m+1}\over dx^{m+1}}(x^{n_1+1}-x^{-n_1-1})|_{x=1},\eqno(4.50)$$
$$\dlt_{r_2,q}((^{-n_2-1}_{\;\;\;\;m})-(^{n_2}_m))={\dlt_{r_2,q}\over m!}{d^m\over dx^m}(x^{-n_2-1}-x^{n_2})|_{x=1}.\eqno(4.51)$$
We can prove the simplicity of $R^{\ast}_{[k_1,k_2]}$ by Theorem 3.2, (4.47)-(4.51) and the same arguments as those in the proof of the simplicity of $R_{[k_1,k_2],\ell+1}$.

The following identities imply the generator property by the arguments in (4.19)-(4.42): by (2.31),
\begin{eqnarray*}& &[(E_{p_1,k_1+p_2})_{[0,1]}+(E_{k_1+p_2,p_1})_{[1,0]}](\es/2)[(E_{q_1,k_1+q_2})_{[0,1]}(m,n)+(E_{k_1+q_2,q_1})_{[1,0]}(n,m)]\\&=&\dlt_{p_1,q_1}(m+1)[(E_{k_1+p_2,k_1+q_2})_{[1,1]}(m+(1-\es)/2,n)\\& &-(E_{k_1+q_2,k_1+p_2})_{[1,1]}(n,m+(1-\es)/2)]+\dlt_{p_2,q_2}[E_{q_1,p_1}(m,n-(\es+1)/2)\\& &+E_{p_1,q_1}(n-(\es+1)/2,m)],\hspace{9.1cm}(4.52)\end{eqnarray*}
\begin{eqnarray*}& &[(E_{p_1,k_1+p_2})_{[0,1]}+(E_{k_1+p_2,p_1})_{[1,0]}](1/2)(E_{p_1,p_1}(m,n)+E_{p_1,p_1}(n,m))\\&=&(n+1)[(E_{p_1,k_1+p_2})_{[0,1]}(m,n)+(E_{k_1+p_2,p_1})_{[1,0]}(n,m)]\\& &+(m+1)[(E_{p_1,k_1+p_2})_{[0,1]}(n,m)+(E_{k_1+p_2,p_1})_{[1,0]}(m,n)],\hspace{4.2cm}(4.53)\end{eqnarray*}
\begin{eqnarray*}& &[(E_{p_1,k_1+p_2})_{[0,1]}+(E_{k_1+p_2,p_1})_{[1,0]}](\es/2)\\& &[(E_{k_1+p_2,k_1+p_2})_{[1,1]}(m,n)-(E_{k_1+p_2,k_1+p_2})_{[1,1]}(n,m)]\\&=&(E_{p_1,k_1+p_2})_{[0,1]}(m+(1-\es)/2,n)+(E_{k_1+p_2,p_1})_{[1,0]}(n,m+(1-\es)/2)\\& &-(E_{p_1,k_1+p_2})_{[0,1]}(n+(1-\es)/2,m)-(E_{k_1+p_2,p_1})_{[1,0]}(m,n+(1-\es)/2),\hspace{1.5cm}(4.54)\end{eqnarray*}
\begin{eqnarray*}& &[(E_{k_1+1,k_1+2})_{[1,1]}(1,0)-(E_{k_1+2,k_1+1})_{[1,1]}(0,1)](0)\\& &[(E_{k_1+1,k_1+2})_{[1,1]}(1,0)-(E_{k_1+2,k_1+1})_{[1,1]}(0,1)]\\&=& 2[(E_{k_1+2,k_1+2})_{[1,1]}(1,0)-(E_{k_1+2,k_1+2})_{[1,1]}(0,1)]\hspace{5.8cm}(4.55)\end{eqnarray*}
if $k_2>1$, for $p_1,q_1\in\ol{1,k_1},\;p_2,q_2\in\ol{1,k_2},\;m,n\in\Bbb{N}$ and $\es=\pm 1.\qquad\Box$
\psp

Assume that
$$k_1=2\ell_1,\;\;\;k_2=2\ell_2\eqno(4.56)$$
 for some positive integers $\ell_1$ and $\ell_2$. Set
$$S=\left(\begin{array}{cccc}& I_{\ell_1}& &\\-I_{\ell_1}&& &\\& & & I_{\ell_2}\\& & -I_{\ell_2}&\end{array}\right),\eqno(4.57)$$
 where the empty entries denote zero matrices. Define a map $\sgm_2: M_{k\times k}(\Bbb{F})\rightarrow M_{k\times k}(\Bbb{F})$ by
$$\sgm_2(A)=S A^tS^{-1}\qquad \for\;\;A\in M_{k\times k}(\Bbb{F}).\eqno(4.58)$$
Then $\sgm_2$ is an involutive anti-isomorphism of $M_{k\times k}(\Bbb{F})$ preserving the $\Bbb{Z}_2$-grading in (4.2) and (4.3). In terms of block matrices, we have:
$$\left(\begin{array}{rrrr}A_{1,1},&A_{1,2},& A_{1,3},&A_{1,4}\\A_{2,1},&A_{2,2},& A_{2,3},&A_{2,4}\\A_{3,1},&A_{3,2},& A_{3,3},&A_{3,4}\\A_{4,1},&A_{4,2},& A_{4,3},&A_{4,4}\end{array}\right)\stl{\sgm_2}{\mapsto}
\left(\begin{array}{rrrr}A_{2,2}^t,&-A_{1,2}^t,& A_{4,2}^t,&-A_{3,2}^t\\-A_{2,1}^t,&A_{1,1}^t,& -A_{4,1}^t,&A_{3,1}^t\\A_{2,4}^t,&-A_{1,4}^t,& A_{4,4}^t,&-A_{3,4}^t\\-A_{2,3}^t,&A_{1,3}^t,&-A_{4,3}^t,&A_{3,3}^t\end{array}\right),\eqno(4.59)$$
where
$$A_{1,1},A_{1,2}, A_{2,1}\in M_{\ell_1\times \ell_1}(\Bbb{F}),\;\; A_{1,3},A_{1,4}\in M_{\ell_1\times \ell_2}(\Bbb{F}),\;\;A_{3,1},A_{4,1}\in M_{\ell_2\times \ell_1}(\Bbb{F}).\eqno(4.60)$$
 Thus we have the following subalgebra of $R_{[k_1, k_2],1}$: $R_{[k_1,k_2]}^{\dg}=(R_{[k_1,k_2]}^{\dg})_0+(R_{[k_1,k_2]}^{\dg})_1$ with
\begin{eqnarray*}& &(R_{[k_1,k_2]}^{\dg})_0=\mbox{span}\:\{E_{p_1,q_1}(m,n)+E_{\ell_1+q_1,\ell_1+p_1}(n,m), E_{p_1,\ell_1+q_1}(m,n)-E_{q_1,\ell_1+p_1}(n,m),\\& &E_{\ell_1+p_1,q_1}(m,n)-E_{\ell_1+q_1,p_1}(m,n),(E_{k_1+p_2,k_1+q_2})_{[1,1]}(m,n)-(E_{k_1+\ell_2+q_2,k_1+\ell_2+p_2})_{[1,1]}(n,m),\\& &(E_{k_1+p_2,k_1+\ell_2+q_2})_{[1,1]}(m,n)+(E_{k_1+q_2,k_1+\ell_2+p_2})_{[1,1]}(n,m),\\& &(E_{k_1+\ell_2+p_2,k_1+q_2})_{[1,1]}(m,n)+(E_{k_1+\ell_2+q_2,k_1+p_2})_{[1,1]}(n,m)\\& &\mid p_1,q_1\in\ol{1,\ell_1},\;p_2,q_2\in\ol{1,\ell_2},\;m,n\in\Bbb{N}\},\hspace{7.1cm}(4.61)\end{eqnarray*}
\begin{eqnarray*}& &(R_{[k_1,k_2]}^{\dg})_1=\mbox{span}\:\{(E_{p_1,k_1+p_2})_{[0,1]}(m,n)+(E_{k_1+\ell_2+p_2,\ell_1+p_1})_{[1,0]}(n,m),\\& &(E_{p_1,k_1+\ell_2+p_2})_{[0,1]}(m,n)-(E_{k_1+p_2,\ell_1+p_1})_{[1,0]}(n,m),\\& &(E_{\ell_1+p_1,k_1+p_2})_{[0,1]}(m,n)-(E_{k_1+\ell_2+p_2,p_1})_{[1,0]}(n,m),\\& &(E_{\ell_1+p_1,k_1+\ell_2+p_2})_{[0,1]}(m,n)+(E_{k_1+p_2,p_1})_{[1,0]}(n,m)\\& &\mid p_1\in\ol{1,\ell_1},\;p_2\in\ol{1,\ell_2},\;m,n\in\Bbb{N}\}.\hspace{8.3cm}(4.62)\end{eqnarray*}

{\bf Theorem 4.3}.  {\it The algebra} $(R^{\dg}_{[k_1,k_2]},\ptl,Y^+(\cdot,z))$ {\it is simple. Moreover, it is generated by} 
\begin{eqnarray*}& &(R_{[k_1,k_2]}^{\dg})^{(3/2)}=\mbox{span}\:\{(E_{p_1,k_1+p_2})_{[0,1]}+(E_{k_1+\ell_2+p_2,\ell_1+p_1})_{[1,0]},\\ & &(E_{p_1,k_1+\ell_2+p_2})_{[0,1]}-(E_{k_1+p_2,\ell_1+p_1})_{[1,0]},(E_{\ell_1+p_1,k_1+p_2})_{[0,1]}-(E_{k_1+\ell_2+p_2,p_1})_{[1,0]},\\& &(E_{\ell_1+p_1,k_1+\ell_2+p_2})_{[0,1]}+(E_{k_1+p_2,p_1})_{[1,0]}\mid p_1\in\ol{1,\ell_1},\;p_2\in\ol{1,\ell_2}\}\hspace{3.8cm}(4.63)\end{eqnarray*}
{\it (cf. (2.61)) if} $k>4$ {\it and by} 
$$(R_{[k_1, k_2]}^{\dg})^{(3/2)}+\Bbb{F}[(E_{3,3})_{[1,1]}(0,1)-(E_{4,4})_{[1,1]}(1,0)]\eqno(4.64)$$
{\it if} $k=4$.
\psp

{\it Proof}. The following eight equations will be used for the proof of the simplicity of $R^{\dg}_{[k_1,k_2]}$: for $r_1,p_1\in\ol{1,\ell_1},\;r_2,p_2\in\ol{1,\ell_2}$ and $m,n_1,n_2\in\Bbb{N}$, by (2.31), we obtain
\begin{eqnarray*}&&(E_{r_1,r_1}(0,m)+E_{\ell_1+r_1,\ell_1+r_1}(m,0))(0)\\& &[
(E_{p_1,k_1+p_2})_{[0,1]}(n_1,n_2)+(E_{k_1+\ell_2+p_2,\ell_1+p_1})_{[1,0]}(n_2,n_1)]\\&=&\dlt_{r_1,p_1}(n_1+1)(^{-n_1-2}_{\;\;\;\;\:m})[(E_{p_1,k_1+p_2})_{[0,1]}(n_1,n_2)+(E_{k_1+\ell_2+p_2,\ell_1+p_1})_{[1,0]}(n_2,n_1)],\hspace{0.9cm}(4.65)\end{eqnarray*}
\begin{eqnarray*}& &(E_{r_1,r_1}(0,m)+E_{\ell_1+r_1,\ell_1+r_1}(m,0))(0)\\& &[(E_{p_1,k_1+\ell_2+p_2})_{[0,1]}(n_1,n_2)-(E_{k_1+p_2,\ell_1+p_1})_{[1,0]}(n_2,n_1)]
\\&=&\dlt_{r_1,p_1}(n_1+1)(^{-n_1-2}_{\;\;\;\;\:m})[(E_{p_1,k_1+\ell_2+p_2})_{[0,1]}(n_1,n_2)-(E_{k_1+p_2,\ell_1+p_1})_{[1,0]}(n_2,n_1)],\hspace{0.9cm}(4.66)\end{eqnarray*}
\begin{eqnarray*}& &(E_{r_1,r_1}(0,m)+E_{\ell_1+r_1,\ell_1+r_1}(m,0))(0)\\& &
[(E_{\ell_1+p_1,k_1+p_2})_{[0,1]}(n_1,n_2)-(E_{k_1+\ell_2+p_2,p_1})_{[1,0]}(n_2,n_1)]\\&=&\dlt_{r_1,p_1}(n_1+1)(^{n_1}_m)[(E_{\ell_1+p_1,k_1+p_2})_{[0,1]}(n_1,n_2)-(E_{k_1+\ell_2+p_2,p_1})_{[1,0]}(n_2,n_1)],\hspace{1.3cm}(4.67)\end{eqnarray*}
\begin{eqnarray*}& &(E_{r_1,r_1}(0,m)+E_{\ell_1+r_1,\ell_1+r_1}(m,0))(0)\\&&
[E_{\ell_1+p_1,k_1+\ell_2+p_2})_{[0,1]}(n_1,n_2)+(E_{k_1+p_2,p_1})_{[1,0]}(n_2,n_1)]\\&=&\dlt_{r_1,p_1}(n_1+1)(^{n_1}_m)[E_{\ell_1+p_1,k_1+\ell_2+p_2})_{[0,1]}(n_1,n_2)+(E_{k_1+p_2,p_1})_{[1,0]}(n_2,n_1)],\hspace{1.6cm}(4.68)\end{eqnarray*}
\begin{eqnarray*}& &[(E_{k_1+r_2,k_1+r_2})_{[1,1]}(0,m)-(E_{k_1+\ell_2+r_2,k_1+\ell_2+r_2})_{[1,1]}(m,0)](0)\\&&[(E_{p_1,k_1+p_2})_{[0,1]}(n_1,n_2)+(E_{k_1+\ell_2+p_2,\ell_1+p_1})_{[1,0]}(n_2,n_1)]\\&=&-\dlt_{r_2,p_2}(^{n_2}_m)
[(E_{p_1,k_1+p_2})_{[0,1]}(n_1,n_2)+(E_{k_1+\ell_2+p_2,\ell_1+p_1})_{[1,0]}(n_2,n_1)],\hspace{2.6cm}(4.69)\end{eqnarray*}
\begin{eqnarray*}& &[(E_{k_1+r_2,k_1+r_2})_{[1,1]}(0,m)-(E_{k_1+\ell_2+r_2,k_1+\ell_2+r_2})_{[1,1]}(m,0)](0)\\&&[(E_{p_1,k_1+\ell_2+p_2})_{[0,1]}(n_1,n_2)-(E_{k_1+p_2,\ell_1+p_1})_{[1,0]}(n_2,n_1)]\\&=&\dlt_{r_2,p_2}(^{-n_2-1}_{\:\;\;\;m})[(E_{p_1,k_1+\ell_2+p_2})_{[0,1]}(n_1,n_2)-(E_{k_1+p_2,\ell_1+p_1})_{[1,0]}(n_2,n_1)],\hspace{2.3cm}(4.70)\end{eqnarray*}
\begin{eqnarray*}& &[(E_{k_1+r_2,k_1+r_2})_{[1,1]}(0,m)-(E_{k_1+\ell_2+r_2,k_1+\ell_2+r_2})_{[1,1]}(m,0)](0)\\&&[(E_{\ell_1+p_1,k_1+p_2})_{[0,1]}(n_1,n_2)-(E_{k_1+\ell_2+p_2,p_1})_{[1,0]}(n_2,n_1)]\\&=&-\dlt_{r_2,p_2}(^{n_2}_m) [(E_{\ell_1+p_1,k_1+p_2})_{[0,1]}(n_1,n_2)-(E_{k_1+\ell_2+p_2,p_1})_{[1,0]}(n_2,n_1)],\hspace{2.6cm}(4.71)\end{eqnarray*}
\begin{eqnarray*}& &[(E_{k_1+r_2,k_1+r_2})_{[1,1]}(0,m)-(E_{k_1+\ell_2+r_2,k_1+\ell_2+r_2})_{[1,1]}(m,0)](0)\\& &[E_{\ell_1+p_1,k_1+\ell_2+p_2})_{[0,1]}(n_1,n_2)+(E_{k_1+p_2,p_1})_{[1,0]}(n_2,n_1)]\\&=&\dlt_{r_2,p_2}(^{-n_2-1}_{\:\;\;\;m})[(E_{\ell_1+p_1,k_1+\ell_2+p_2})_{[0,1]}(n_1,n_2)+(E_{k_1+p_2,p_1})_{[1,0]}(n_2,n_1)].\hspace{2.5cm}(4.72)\end{eqnarray*}

The following fourteen equations will be used both for the proof of the simplicity and generator property of $R^{\dg}_{[k_1,k_2]}$: for $p_1,q_1\in\ol{1,\ell_1},\;p_2,q_2\in\ol{1,\ell_2}$ and $n_1,n_2\in\Bbb{N}$, by (2.31), (2.61) and (2.62), we obtain
\begin{eqnarray*}&&[(E_{p_1,k_1+p_2})_{[0,1]}+(E_{k_1+\ell_2+p_2,\ell_1+p_1})_{[1,0]}](1/2)(E_{p_1,p_1}(n_1,n_2)+E_{\ell_1+p_1,\ell_1+p_1}(n_2,n_1))\\&=&(n_2+1)[(E_{p_1,k_1+p_2})_{[0,1]}(n_1,n_2)+(E_{k_1+\ell_2+p_2,\ell_1+p_1})_{[1,0]}(n_2,n_1)],\hspace{3.1cm}(4.73)\end{eqnarray*}
\begin{eqnarray*}& &[(E_{p_1,k_1+\ell_2+p_2})_{[0,1]}-(E_{k_1+p_2,\ell_1+p_1})_{[1,0]}](1/2)(E_{p_1,p_1}(n_1,n_2)+E_{\ell_1+p_1,\ell_1+p_1}(n_2,n_1))
\\&=&(n_2+1)[(E_{p_1,k_1+\ell_2+p_2})_{[0,1]}(n_1,n_2)-(E_{k_1+p_2,\ell_1+p_1})_{[1,0]}(n_2,n_1)],\hspace{3.1cm}(4.74)\end{eqnarray*}
\begin{eqnarray*}& &[(E_{\ell_1+p_1,k_1+p_2})_{[0,1]}-(E_{k_1+\ell_2+p_2,p_1})_{[1,0]}](1/2)(E_{p_1,p_1}(n_1,n_2)+E_{\ell_1+p_1,\ell_1+p_1}(n_2,n_1))\\&=&(n_1+1)[(E_{\ell_1+p_1,k_1+p_2})_{[0,1]}(n_2,n_1)-(E_{k_1+\ell_2+p_2,p_1})_{[1,0]}(n_1,n_2)],\hspace{3cm}(4.75)\end{eqnarray*}
\begin{eqnarray*}& &[(E_{\ell_1+p_1,k_1+\ell_2+p_2})_{[0,1]}+(E_{k_1+p_2,p_1})_{[1,0]}](1/2)(E_{p_1,p_1}(n_1,n_2)+E_{\ell_1+p_1,\ell_1+p_1}(n_2,n_1))
\\&=&(n_1+1)[(E_{\ell_1+p_1,k_1+\ell_2+p_2})_{[0,1]}(n_2,n_1)+(E_{k_1+p_2,p_1})_{[1,0]}(n_1,n_2)],\hspace{3.1cm}(4.76)\end{eqnarray*}
\begin{eqnarray*}& &[(E_{p_1,k_1+p_2})_{[0,1]}+(E_{k_1+\ell_2+p_2,\ell_1+p_1})_{[1,0]}](-1/2)\\& &[(E_{k_1+p_2,k_1+p_2})_{[1,1]}(n_1,n_2)-(E_{k_1+\ell_2+p_2,k_1+\ell_2+p_2})_{[1,1]}(n_2,n_1)]\\&=&(E_{p_1,k_1+p_2})_{[0,1]}(n_1,n_2)+(E_{k_1+\ell_2+p_2,\ell_1+p_1})_{[1,0]}(n_2,n_1),\hspace{4.7cm}(4.77)\end{eqnarray*}
\begin{eqnarray*}& &[(E_{p_1,k_1+\ell_2+p_2})_{[0,1]}-(E_{k_1+p_2,\ell_1+p_1})_{[1,0]}](-1/2)\\& &[(E_{k_1+p_2,k_1+p_2})_{[1,1]}(n_1,n_2)-(E_{k_1+\ell_2+p_2,k_1+\ell_2+p_2})_{[1,1]}(n_2,n_1)]\\&=&-[(E_{p_1,k_1+\ell_2+p_2})_{[0,1]}(n_2,n_1)-(E_{k_1+p_2,\ell_1+p_1})_{[1,0]}(n_1,n_2)],\hspace{4.1cm}(4.78)\end{eqnarray*}
\begin{eqnarray*}& &[(E_{\ell_1+p_1,k_1+p_2})_{[0,1]}-(E_{k_1+\ell_2+p_2,p_1})_{[1,0]})](-1/2)\\& &[(E_{k_1+p_2,k_1+p_2})_{[1,1]}(n_1,n_2)-(E_{k_1+\ell_2+p_2,k_1+\ell_2+p_2})_{[1,1]}(n_2,n_1)]\\&=&(E_{\ell_1+p_1,k_1+p_2})_{[0,1]}(n_1,n_2)-(E_{k_1+\ell_2+p_2,p_1})_{[1,0]}(n_2,n_1),\hspace{4.7cm}(4.79)\end{eqnarray*}
\begin{eqnarray*}& &[(E_{\ell_1+p_1,k_1+\ell_2+p_2})_{[0,1]}+(E_{k_1+p_2,p_1})_{[1,0]}](-1/2)\\& &[(E_{k_1+p_2,k_1+p_2})_{[1,1]}(n_1,n_2)-(E_{k_1+\ell_2+p_2,k_1+\ell_2+p_2})_{[1,1]}(n_2,n_1)]\\&=&-[(E_{\ell_1+p_1,k_1+\ell_2+p_2})_{[0,1]}(n_2,n_1)+(E_{k_1+p_2,p_1})_{[1,0]}(n_1,n_2)],\hspace{4.2cm}(4.80)\end{eqnarray*}
\begin{eqnarray*}& &[(E_{p_1,k_1+p_2})_{[0,1]}+(E_{k_1+\ell_2+p_2,\ell_1+p_1})_{[1,0]}](\es/2)\\& &[(E_{\ell_1+q_1,k_1+\ell_2+q_2})_{[0,1]}(n_1,n_2)+(E_{k_1+q_2,q_1})_{[1,0]}(n_2,n_1)]\\&=&\dlt_{p_2,q_2}(E_{p_1,q_1}(n_2-(\es+1)/2,n_1)+
E_{\ell_1+q_1,\ell_1+p_1}(n_1,n_2-(\es+1)/2))\\& &-\dlt_{p_1,q_1}(n_1+1)[(E_{k_1+q_2,k_1+p_2})_{[1,1]}(n_2,n_1+(1-\es)/2)\\& &-(E_{k_1+\ell_2+p_2,k_1+\ell_2+q_2})_{[1,1]}(n_1+(1-\es)/2,n_2)],\hspace{5.9cm}(4.81)\end{eqnarray*}
\begin{eqnarray*}& &[(E_{p_1,k_1+p_2})_{[0,1]}+(E_{k_1+\ell_2+p_2,\ell_1+p_1})_{[1,0]}](-1/2)\\& &[(E_{q_1,k_1+\ell_2+q_2})_{[0,1]}(n_1,n_2)-(E_{k_1+q_2,\ell_1+q_1})_{[1,0]}(n_2,n_1)]\\&=&
-\dlt_{p_2,q_2}(E_{p_1,\ell_1+q_1}(n_2,n_1)-
E_{\ell_1+q_1,\ell_1+p_1}(n_1,n_2)),\hspace{5.5cm}(4.82)\end{eqnarray*}
\begin{eqnarray*}& &[(E_{p_1,k_1+p_2})_{[0,1]}+(E_{k_1+\ell_2+p_2,\ell_1+p_1})_{[1,0]}](1/2)\\& &[(E_{\ell_1+q_1,k_1+q_2})_{[0,1]}(n_1,n_2)-(E_{k_1+\ell_2+q_2,q_1})_{[1,0]})(n_2,n_1)]\\&=&\dlt_{p_1,q_1}(n_1+1)[(E_{k_1+\ell_2+p_2,k_1+q_2})_{[1,1]}(n_1,n_2)+(E_{k_1+\ell_2+q_2,k_1+p_2})_{[1,1]}(n_2,n_1)],\hspace{1.1cm}(4.83)\end{eqnarray*}
\begin{eqnarray*}& &[E_{\ell_1+p_1,k_1+\ell_2+p_2})_{[0,1]}+(E_{k_1+p_2,p_1})_{[1,0]}](\es/2)\\& &
[(E_{q_1,k_1+q_2})_{[0,1]}(n_1,n_2)+(E_{k_1+\ell_2+q_2,\ell_1+q_1})_{[1,0]}(n_2,n_1)]\\&=&\dlt_{p_2,q_2}(E_{q_1,p_1}(n_1,n_2-(\es+1)/2)+
E_{\ell_1+p_1,\ell_1+q_1}(n_2-(\es+1)/2),n_1)\\& &+\dlt_{p_1,q_1}(n_1+1)[(E_{k_1+p_2,k_1+q_2})_{[1,1]}(n_1+(1-\es)/2,n_2)\\& &-(E_{k_1+\ell_2+q_2,k_1+\ell_2+p_2})_{[1,1]}(n_2,n_1+(1-\es)/2)],\hspace{5.8cm}(4.84)\end{eqnarray*}
\begin{eqnarray*}& &[E_{\ell_1+p_1,k_1+\ell_2+p_2})_{[0,1]}+(E_{k_1+p_2,p_1})_{[1,0]}](1/2)\\& &[(E_{q_1,k_1+\ell_2+q_2})_{[0,1]}(n_1,n_2)-(E_{k_1+q_2,\ell_1+q_1})_{[1,0]}(n_2,n_1)]\\&=&
\dlt_{p_1,q_1}(n_1+1)[(E_{k_1+p_2,k_1+\ell_2+q_2})_{[1,1]}(n_1,n_2)-(
E_{k_1+q_2,k_1+\ell_2+p_2})_{[1,1]}(n_2,n_1)],\hspace{1.1cm}(4.85)\end{eqnarray*}
\begin{eqnarray*}& &[E_{\ell_1+p_1,k_1+\ell_2+p_2})_{[0,1]}+(E_{k_1+p_2,p_1})_{[1,0]}](-1/2)\\& &[(E_{\ell_1+q_1,k_1+q_2})_{[0,1]}(n_1,n_2)-(E_{k_1+\ell_2+q_2,q_1})_{[1,0]})(n_2,n_1)]\\&=&\dlt_{p_2,q_2}[E_{\ell_1+q_1,p_1}(n_1,n_2)-E_{\ell_1+p_1,q_1}(n_2,n_1)].\hspace{6.5cm}(4.86)\end{eqnarray*}

If $\ell_2>1$, we need the following equation for the proof of the generator property of $R^{\dg}_{[k_1,k_2]}$:
\begin{eqnarray*}& &[(E_{k_1+1,k_1+2})_{[1,1]}(0,1)-(E_{k_1+\ell_2+2,k_1+\ell_2+1})_{[1,1]}(1,0)](0)\\& &[(E_{k_1+2,k_1+1})_{[1,1]}(0,1)-(E_{k_1+\ell_2+1,k_1+\ell_2+2})_{[1,1]}(1,0)]\\&=&-(E_{k_1+1,k_1+1})_{[1,1]}(0,1)+(E_{k_1+\ell_2+1,k_1+\ell_2+1})_{[1,1]}(1,0)\\& &-(E_{k_1+2,k_1+2})_{[1,1]}(0,1)+(E_{k_1+\ell_2+2,k_1+\ell_2+2})_{[1,1]}(1,0)\hspace{4.9cm}(4.87)\end{eqnarray*}
by (2.31) and (2.62).

Now the simplicity of $R^{\dg}_{[k_1,k_2]}$ follows from Theorem 3.3, (4.65)-(4.84) and the same arguments as those in the proof of the simplicity of $R_{[k_1,k_2],\ell+1}$. Moreover, the generator property can be obtained by (4.73)-(4.87) and similar arguments as those in (4.19)-(4.42).$\qquad\Box$

\newpage

\noindent{\Large \bf References}

\hspace{0.5cm}

\begin{description}

\item[{[BVV]}] B. Bakalov, V. G. Kac and A. Voronov, Cohomology of conformal algebras, {\it Commun. Math. Phys.} {\bf 200} (1999), 561-598.

\item[{[Ba]}] I. Bakas, The large-$N$ limit of extended conformal symmetries, {\it Phys. Lett. B.} {\bf 228} (1989), 57-63.

\item[{[BKV]}] B. Bakalov, V. G. Kac and A. A. Voronov, Cohomology of conformal algebras, {\it Commun. Math. Phys.} {\bf 200} (1999), 561-598.

\item[{[CK1]}] S.-J. Cheng and V. G. Kac, A new $N=6$ superconformal algebras, {\it Commun. Math. Phys.} {\bf 186} (1997), 219-231.

\item[{[CK2]}] ---, Conformal modules, {\it Asian J. Math.} {\bf 1} (1997), 181-193.

\item[{[CK3]}] ---, Erratum: ``Conformal modules,'' {\it Asian J. Math.} {\bf 2} (1998), 153-156.

\item[{[CKW]}] S.-J. Cheng, V. G. Kac and M. Wakimoto, Extensions of conformal modules, {\it Topological field theory, primitive forms and related topics (Kyoto, 1996)}, 79-129. Progr. Math., {\bf 160}, Birkh\"{a}user Boston, Boston, MA, 1998. 

\item[{[DK]}] A. D'Andrea and V. G. Kac, Structure theory of finite conformal algebras, {\it Selecta Math (N.S.)} {\bf 4} (1998), 377-418.

\item[{[DW]}] A. Das and  W.-J. Huang, The Hamiltonian structures associatted with a generalized Lax operators, {\it J. Math. Phys.} {\bf 33} (7) (1992), 2487-2497.   

\item[{[DHP]}] A. Das, W.-J. Huang and S. Panda, The Hamiltonian structures of the KP Hierarchy, {\it Phys. Lett.} {\bf B271} (1991), 109-115.

\item[{[FKW]}] E. Frenkel, V. Kac and W. Wang, ${\cal W}_{1+\infty}$ and ${\cal W}(gl_N)$ with central charge $N$, {\it Commun. Math. Phys.} {\bf 170} (1995), 337-357.

\item[{[GK]}] M. I. Golenishcheva-Kutuzova and V. G. Kac, $\G$-conformal algebras, {\it J. Math. Phys.} {\bf 39} (1998), no. 4, 2290-2305.

\item[{[K1]}] V. G. Kac, Simple graded Lie algebras of finite growth,  {\it Funct. Anal. Appl.} {\bf 1} (1967), 328-329.

\item[{[K2]}] ---, Lie superalgebras, {\it Adv. Math.} {\bf 26} (1977), 8-96.

\item[{[K3]}] ---, {\it Vertex algebras for beginners}, University lectures series, Vol {\bf 10}, AMS. Providence RI, 1996.

\item[{[K4]}] ---, Superconformal algebras and transitive group actions on quadrics, {\it Commun. Math. Phys.} {\bf 186} (1997), 233-252.

\item[{[K5]}] ---, Idea of locality, {\it Physical Applications and Mathematical Aspects of Geometry, Groups and Algebras}, Doebener et al eds., World Scientific Publishers, 1997, 16-32.

\item[{[KL]}]  V. G. Kac and J. W. Leur, On classification of superconformal algebras, in S. J. et al. eds. {\it String 88}, World Sci. (1989), 77-106.

\item[{[KT]}] V. G. Kac and I. T. Todorov, Superconformal current algebras and their unitary representations, {\it Commun. Math. Phys}. {\bf 102} (1985), 337-347.

\item[{[KWY]}] V. G. Kac, W. Wang and C. H. Yan, Quasifinite representations of classical Lie subalgebras of ${\cal W}_{1+\infty}$, {\it Adv. Math.} {\bf 139} (1998), 56-140. 

\item[{[V]}] J. W. van de Leur, the $W_{1+\infty}(gl_2)$-symmetries of the $s$-component KP hierarchy, {\it J. Math. Phys.} {\bf 37} (1996), 2315-2337.

\item[{[X1]}] X. Xu, Hamiltonian operators and associative algebras with a derivation, {\it  Lett. Math. Phys.} {\bf 33} (1995), 1-6.

\item[{[X2]}] ---, Hamiltonian superoperators, {\it J. Phys A: Math. \& Gen.} {\bf 28} No. 6 (1995).

\item[{[X3]}] ---, Analogue of the vertex operator triality for ternary codes, {\it J. Pure and Appl. Algebra} {\bf 99} (1995), 53-111 (preprint was circulated in 1991)

\item[{[X4]}] ---, Variational calculus of supervariables and related algebraic structures, {\it J. Algebra}, in press; preprint was circulated in January 1995.

\item[{[X5]}] ---, {\it Introduction to Vertex Operator Superalgebras and Their Modules}, Kluwer Academic Publishers, Dordrecht/Boston/London, 1998.

\item[{[X6]}] ---, Quadratic conformal superalgebras, {\it J. Algebra}, to appear.

\item[{[Y]}] K. Yamagishi, A Hamiltonian structure of $KP$ hierarchy, $W_{1+\infty}$ algebra, and self-dual gravity, {\it Phys. Lett.} {\bf B 259} (1991), 436-441.

\item[{[YW]}] F. Yu and Y.-S. Wu, Hamiltonian structure (anti-) self-adjoint flows in the KP hierarchy and the $W_{1+\infty}$ and $W_{\infty}$ algebras, {\it Phys. Lett.} {\bf B263} (1991), 220-226.

\end{description}
\end{document}